\documentclass[11pt]{nyjm}
\usepackage{mathtools, nccmath}
\usepackage{latexsym}
\usepackage{amsfonts}
\usepackage{amssymb}
\usepackage{psfrag}
\usepackage{graphicx}
\usepackage{epsfig}
\usepackage{amsmath,amsfonts,amsthm}
\usepackage{mathrsfs}
\usepackage{enumerate}
\usepackage{dsfont}
\usepackage{hyperref}
\usepackage{pst-node}
\usepackage{cleveref}
\usepackage{tikz-cd}
\usepackage{multirow}
\usepackage{tikz}
\usetikzlibrary{arrows,decorations.markings}
\usepackage{bookmark}
\usepackage{hyperref}
\hypersetup{
     colorlinks   = true,
     citecolor    = blue
}

\newcommand\restr[2]{{
  \left.\kern-\nulldelimiterspace 
  #1 
  \littletaller 
  \right|_{#2} 
  }}

\newcommand{\littletaller}{\mathchoice{\vphantom{\big|}}{}{}{}}

\newcommand{\f}{\mathcal{F}}

\newcommand{\newsection}[1]{\setcounter{equation}{0}
\setcounter{dfn}{0}
\section{#1}}

\renewcommand{\theequation}{\thesection.\arabic{equation}}

\newtheorem{dfn}{Definition}[section]
\newtheorem{thm}[dfn]{Theorem}
\newtheorem{lmma}[dfn]{Lemma}
\newtheorem{ppsn}[dfn]{Proposition}
\newtheorem{crlre}[dfn]{Corollary}

\newtheorem{rmrk}[dfn]{Remark}

\newcommand{\bbc}{\mathbb{C}}
\newcommand{\bbz}{\mathbb{Z}}
\newcommand{\tr}{\mathrm{tr}}

\def \qed { \mbox{}\hfill
$\Box$\vspace{1ex}}

\title{Sections and Chapters}

\begin{document}
%%%%%%%%%%%%%%%%%%%%%%%%%%%%%%%%%
%%%%%%%%%%%%%%%%%%%%%%%%%%%%%%%%%

\author{\sc{Keshab Chandra Bakshi, Satyajit Guin, Sruthymurali}}
\title{Fourier-theoretic inequalities for inclusions of simple $C^*$-algebras}
\maketitle

%%%%%%%%%%%%%%%%%%%%%%%%%%%%%%%%%%
%%%%%%%%%%%   ABSTRACT    %%%%%%%%%%%%%%%%
%%%%%%%%%%%%%%%%%%%%%%%%%%%%%%%%%%

\begin{abstract}
This paper originates from a naive attempt to establish various non-commutative Fourier-theoretic inequalities  for an inclusion of simple $C^*$-algebras equipped with a conditional expectation of index-finite type. In this setting, we discuss the Hausdorff-Young inequality and Young's inequality. As a consequence, we  prove the Hirschman–Beckner uncertainty principle and Donoho–Stark uncertainty principle. Our results generalize some of the results of Jiang, Liu and Wu [Noncommutative uncertainty principle, J. Funct. Anal., 270(1): 264--311, 2016].
\end{abstract}
\bigskip

{\bf AMS Subject Classification No.:} {\large 46}L{\large 37}\,, {\large 47}L{\large 40}\,, {\large 46}L{\large 05}\,, {\large 43}A{\large 30}\,.

\bigskip
\hypersetup{linkcolor=black}
\tableofcontents

%%%%%%%%%%%%%%%%%%%%%%%%%%%%%%%%%%
%%%%%%%%%   INTRODUCTION    %%%%%%%%%%%%%%%
%%%%%%%%%%%%%%%%%%%%%%%%%%%%%%%%%%

\newsection{Introduction}

Jones had discovered a notion of index $[M:N]$ for an inclusion of type $II_1$ factors $N\subset M$ in \cite{Jo}, which is now an active area of research in operator algebra having applications in various other fields of mathematics and mathematical physics. Subsequently, Kosaki \cite{kosaki} introduced a notion of index for an inclusion of type $III$ factors.  Both type $II_1$ and ($\sigma$-finite) type $III$ factors are particular cases of more general objects, called simple $C^*$-algebras (i.e., $C^*$-algebras having no proper closed ideals). Thus, inclusion of simple $C^*$-algebras encompasses both type $II_1$ and type $III$ subfactor theory. As a generalization of both Jones' and Kosaki's indices, Watatani \cite{Watataniindex} discussed index for an inclusion of $C^*$-algebras with a conditional expectation having a `quasi-basis', a generalization of Pimsner-Popa basis \cite{PP}.  Given a subfactor $N\subset M$ with $[M:N]<\infty,$ Jones crucially observed that one can obtain another type $II_1$ factor $M_1$, called the basic construction, so that $[M_1:M]=[M:N]$  and furthermore, this operation can be iterated to obtain a tower of basic constructions: $N\subset M\subset M_1\subset\cdots\subset M_k\subset \cdots$. This was the key observation in establishing the famous `Jones' index rigidity' result in \cite{Jo}. It is well-known from the early days of subfactor theory that the higher relative commutants $N^{\prime}\cap M_k$ and $M^{\prime}\cap M_k$ have incredibly  rich structures. Using the relative commutants, Popa had associated a `standard invariant' to a subfactor~: the $\lambda$-lattice \cite{P}, which is arguably the most powerful invariant of a subfactor. Subsequently, Jones discovered a pictorial description of the standard invariant what he called `planar algebra' \cite{Jo2} and it becomes an indispensable tool in subfactor theory. In another direction, Ocneanu introduced a (fundamental) notion of Fourier transform $\mathcal{F}$ from $N^{\prime}\cap M_1$ onto $M^{\prime}\cap M_2$ (see also \cite{Bi94, BiJo2} for details). This  generalizes the classical notion of Fourier transform for a finite abelian group. Furthermore, using the Fourier transform, one can associate a new multiplication structure on $N^{\prime}\cap M_1$ that generalizes the classical convolution. Fourier transform on the relative commutants plays a major role in the abstract subfactor theory; for example, the Fourier transform and `rotation operators' are instrumental in the formalism of Jones' planar algebra, Ocneanu's paragroups and Popa's $\lambda$-lattice. Fourier transform also appeared naturally in Bisch's biprojection theory (see \cite{Bisch}) which is an indispensable tool in the theory of intermediate subfactors. In Jones' planar algebraic language, the Fourier transform, convolution and rotation operator have beautiful pictorial descriptions (see \cite{Jo2,BiJo2}). Exploiting these pictorial formulations, in a recent paper \cite{Liunoncommutative} Jiang, Liu and Wu provided a non-commutative version of the Hausdorff-Young inequality, the Young's inequality and uncertainty principles for a subfactor $N\subset M$ with $[M:N]<\infty$ and $N^{\prime}\cap M=\mathbb{C}$ (such a subfactor is called irreducible). Moreover, for any extremal subfactor which is not necessarily irreducible, these inequalities were proved in the Section 7 of \cite{liu2016exchange} using planar algebraic techniques. However, the proofs use sphericality of the planar algebra and so are no longer valid for the non-extremal subfactors. We also mention some related works for Kac algebras and locally compact quantum groups, see for instance \cite{liuquantumgroup,liukacalgebra}.
\smallskip

Analogous to Jones' subfactor theory, given a unital inclusion of simple $C^*$-algebras $B\subset A$, recently in \cite{BakshiVedlattice}, one of the authors and Gupta have systematically developed a Fourier theory on the relative commutants $B^{\prime}\cap A_k$ and $A^{\prime}\cap A_k$, where $B\subset A\subset A_1\subset A_2\subset\cdots\subset A_k\subset \cdots$ is Watatani's tower of $C^*$-basic constructions (this notion parallels to the Jones' tower of basic construction in the $C^*$-world). More precisely, we have a notion of the Fourier transform $\mathcal{F}:B^{\prime}\cap A_1\to A^{\prime}\cap A_2$, rotation map ${\rho}_{+}:B^{\prime}\cap A_1\to B^{\prime}\cap A_1$ and the convolution product $*:B^{\prime}\cap A_1\to B^{\prime}\cap A_1$. The crucial ingredient that was used repeatedly is that of minimal conditional expectation and the minimal (Watatani) index ${[A:B]}_0$. In particular, for a subfactor $N\subset M$ with $[M:N]<\infty$ we have another notion of index ${[M:N]}_0$ and it is a fact that $N\subset M$ is `extremal' if and only if $[M:N]={[M:N]}_0$. In this paper, we prove a non-commutative version of the Hausdorff-Young inequality, Young's inequality and a couple of uncertainty principles in the setting of inclusion of simple $C^*$-algebras and in particular, for a subfactor (of both type $II_1$ factors and ($\sigma$-finite) type $III$-factors) with finite index that is not-necessarily extremal. Unlike \cite{Liunoncommutative}, in the $C^*$-setting the main difficulty lies in the fact that it is still not known whether there are pictorial descriptions (similar to planar algebra) of the Fourier transform and convolution. However, we have found a way around using the relationship between quasi-basis and minimal conditional expectations as exploited in \cite{BakshiVedlattice}. Our proofs are inspired by the corresponding proofs in \cite{Liunoncommutative}.

%%%%%%%%%%%%%%%%%%%%%%%%%%%%%%%
%%%%%%%            Section 2            %%%%%%%%%%%%%
%%%%%%%%%%%%%%%%%%%%%%%%%%%%%%%

\newsection{Preliminaries}\label{Sec 2}

In this section, we fix the notations and briefly recall a few key ingredients which we repeatedly use in the sequel. For more details the readers are requested to see \cite{Watataniindex, KajiwaraWatatani, BakshiVedlattice}.

\subsection{Non-commutative conditional expectation}

Analogous to the classical case, in the theory of operator algebras we have a well-studied notion of conditional expectation. Suppose we have an inclusion of $C^*$-algebras $B\subset A$. All $C^*$-algebras considered in this article will be unital, and all inclusion $B\subset A$ of $C^*$-algebras will be considered as unital inclusion. A conditional expectation $E: A\to B$ is a linear surjective map satisfying
\begin{center}
$E(ba)=bE(a)\quad,\quad E(ab)=E(a)b\quad\text{and}\quad E(b)=b$
\end{center}
for all $b\in B$ and $a\in A$. In particular, $E$ is a norm one projection (see for instance in \cite{kadisonschwarz}). A $C^*$-inclusion may not have any conditional expectation. However, if we consider an inclusion of von Neumann algebras $N\subset M$ with $M$ having a tracial state $\tr$ (i.e., a $\sigma$- weak-operator-topology(WOT) continuous linear functional $\tr:M\to\mathbb{C}$ satisfying $\tr(xy)=\tr(yx)$ for all $x,y\in M$ and $\tr(1)=1$), there always exists a unique `trace preserving' conditional expectation, denoted by $E^M_N$. More precisely, $E^M_N$ is characterized by $\tr(nE^M_N(m))=\tr(nm)$ for all $n\in N$ and $m\in M$.

Let us recall the very useful Kadison-Schwarz inequality involving conditional expectation which we shall use later.

\begin{lmma}[\cite{kadisonschwarz}]\label{KSinequality}
Suppose that $N$ and $M$ are von Neumann algebras acting on a Hilbert space $\mathcal{H}$, $N \subset  M$, and $E:M\to N$ is a conditional expectation. Then, for all $x \in M$ one has $E(x)^*E(x)\leq E(x^*x)$.
\end{lmma}

We remark that all conditional expectations in this paper are assumed to be faithful.

\subsection{A quick look at $C^*$-index theory}\label{wat}

Motivated by the Jones' index theory, Watatani developed a theory of index for inclusion of $C^*$-algebras. Given a pair $B \subset A$ of $C^*$-algebras, a conditional expectation $E:A\to B$ is said to be of \textit{index-finite type} if there exists a finite set $\{\lambda_1,\ldots,\lambda_n\}\subset A$ such that
\begin{center}
$x=\sum_{i=1}^nE(x\lambda_i)\lambda^*_i=\sum_{i=1}^n\lambda_iE(\lambda^*_ix)$
\end{center}
for every $x\in A$. The set $\{\lambda_1,\ldots,\lambda_n\}$ is called a \textit{quasi-basis} for $E$ (see \cite{Watataniindex}). The Watatani index of $E$ is defined by
\begin{center}
$\mathrm{Ind}_w (E)= \sum_{i=1}^n \lambda_i\lambda^*_i,$
\end{center}
and is independent of a quasi-basis. In general, $\mathrm{Ind}_w(E)$ is not a scalar but an invertible positive element in  the center $\mathcal{Z}(A)$ of $A$. In particular, if $A$ is a simple $C^*$-algebra, the index is scalar-valued. We denote by $\mathcal{E}_0(A,B)$ the set of all index-finite type conditional expectations from $A$ onto $B$. A conditional expectation $F\in\mathcal{E}_0( A , B)$ is said to be minimal if it satisfies $\mathrm{Ind}_w(F) \leq \mathrm{Ind}_w(E)$ for all $E\in\mathcal{E}_0(A, B)$ (see \cite{Watataniindex} and the references therein). For inclusion of simple $C^*$-algebras, we have a privileged minimal conditional expectation as mentioned below.

\begin{thm}{\em \cite[Theorem $2.12.3$]{Watataniindex}} Let $B\subset A$ be an inclusion of simple $C^*$-algebras such that $\,\mathcal{E}_0(A,B)\neq\emptyset$. Then, there exists a unique minimal conditional expectation $E_0$ from $A$ onto $B$.
\end{thm}

For inclusion of simple $C^*$-algebras $B\subset A$, the minimal index is defined as
\begin{center}
$[A:B]_0:=\mathrm{Ind}_w(E_0).$
\end{center}
As is customary, we shall denote ${[A:B]}_0$ by $\delta^2$. We point out here that if $N\subset M$ is a subfactor with finite Jones index $[M:N]$ and is irreducible (i.e., $N^{\prime}\cap M=\mathbb{C}$), then the trace preserving conditional expectation $E^M_N$ is the minimal conditional expectation with $[M:N]={[M:N]}_0$. In general, the minimal index and Jones index need not coincide. Indeed, a subfactor is extremal if and only if $[M:N]={[M:N]}_0$. We also remark that irreducibility of a subfactor automatically implies  extremality.

In the subfactor theory, Jones' basic construction plays a pivotal role. Using the language of the Hilbert $C^*$-module, Watatani proposed a parallel notion of basic construction in the $C^*$-world, the so-called $C^*$-basic construction. For the convenience of the reader we briefly recall it here and the details can be found in \cite{Watataniindex}. Let $B \subset A$ be an inclusion of  $C^*$-algebras and $E_B \in \mathcal{E}_0(A,B)$. Then, $A$ is a Hilbert $B$-module with respect to the $B$-valued inner product given by
\begin{equation}\label{B-valued}
\langle x, y\rangle_{B}=E_B(x^*y)\quad \text{ for all } x, y\in A.
\end{equation}
Recall that the space $\mathcal{L}_{B}({A})$ consisting of adjointable $B$-linear maps on ${A}$ is a $C^*$-algebra. For each $a \in A$, consider $\lambda(a)\in\mathcal{L}_B({A})$ given by $\lambda(a)(x)=ax$ for $x\in A$. For $x\in A$, the association $x\mapsto E_B(x)$ is an adjointable projection on ${A}$, and is denoted by $e_B\in \mathcal{L}_B({A})$. The projection $e_B$ is called the Jones projection for the pair $B\subset A$. The $C^*$-basic construction $C^*\langle A,e_B\rangle$ is defined to be the $C^*$-subalgebra generated by $\{\lambda(A),e_B\}$ in $\mathcal{L}_B({A})$. It turns out that $C^*\langle A, e_B\rangle$ equals the closure of the linear span of $\{\lambda(x)e_B\lambda(y):x,y\in A\}$ in the $C^*$-algebra $\mathcal{L}_B({A})$; $\lambda$ is an injective $*$-homomorphism, and thus we may consider $A$ as a $C^*$-subalgebra of $C^*\langle A,e_B\rangle$. It is customary to denote the basic construction by $A_1$ if the conditional expectation is understood from the context, and it is a fact that $A_1$ is  simple whenever $B\subset A$ is an inclusion of simple $C^*$-algebras. There exists a unique finite index conditional expectation $\widetilde{E}_B: A_1\to A$ called the `dual conditional expectation' of $E_B$. Furthermore, $\mathrm{Ind}_w (E_B)=\mathrm{Ind}_w(\widetilde{E}_B)$.
 
\subsection{Tracial states on the relative commutants}\label{tracialstate}
 
Let $B\subset A$ be an inclusion of simple $C^*$-algebras and $E\in\mathcal{E}_0(A,B)$. Let $E_0$ be the unique minimal conditional expectation from $A$ onto $B$. Suppose $B\subset A\subset A_1$ is the basic construction corresponding to $E_0$. Note that ${[A_1:A]}_0={[A:B]}_0$. The dual conditional expectation $\widetilde{E}_0$ is also minimal \cite{Watataniindex}. We put $E_1=\widetilde{E}_0$. Iterating the tower of $C^*$-basic construction for the inclusion $B\subset A$, we obtain
\begin{center}
$B \subset A \subset A_1 \subset A_2 \subset \cdots \subset A_k \subset \cdots$
\end{center}
with unique (dual) minimal conditional expectations $E_{k}: A_{k}\to A_{k-1}$, $k \geq 0$, with the convention that $A_{-1}:=B$ and $A_0:=A$. For each $k\geq 0$, let $e_k$ be the Jones projection in $A_k$. The following extremely useful lemma is the $C^*$-analogue of the `push-down lemma' in subfactor theory \cite{PP}.

\begin{lmma}[\cite{BakshiVedlattice}]\label{pushdown}
If $x_1\in A_1$, then there exists a unique $x_0\in A$ such that $x_1e_1=x_0e_1$, where $x_0=[A:B]_0 E_1(x_1e_1)$.
\end{lmma}

Let $B^\prime\cap A_k:= \{x\in A_k: xb=bx\mbox{ for all }b\in B\}$ be the relative commutants of $B$ in $A_k$. It is known that for each $k$, the relative commutants are finite dimensional \cite{Watataniindex}. On each $B'\cap A_k$, using the minimal conditional expectations, one obtains a consistent `Markov type trace' (Proposition $2.21$ in \cite{BakshiVedlattice}). More precisely, for each $k\geq 0$, the map $\mathrm{tr}_k: B^{\prime}\cap A_k\to\mathbb{C}$ defined by $\mathrm{tr}_k=(E_0\circ E_1\circ\cdots\circ E_k)_{|_{B^{\prime}\cap A_k}}$ becomes a faithful traical state on $B^{\prime}\cap A_k$. 

\begin{ppsn}[\cite{BakshiVedlattice}]\label{m1}
For each $k \geq 0$, $B^{\prime}\cap A_k$ admits a faithful tracial state $\mathrm{tr}_k$ such that
\begin{equation}
\mathrm{tr}_k(xe_k)={\delta}^{-2}\mathrm{tr}_{k-1}(x) \,\,\text{ for all } x\in B^{\prime}\cap A_{k-1}
\end{equation}
and $\restr{\tr_k}{B'\cap A_{k-1}}=\tr_{k-1}$ for all $k\geq 1$.
\end{ppsn}

We shall sometimes drop $k$ and denote $\tr_k$ simply by $\tr$ for notational brevity. The following lemma is very useful.

\begin{lmma}[\cite{BakshiVedlattice}]\label{f2} Let $\{\lambda_i:1\leq i\leq n\}\subset A$ be a quasi-basis for the minimal conditional expectation $E_0$. Then, the $\mathrm{tr}$-preserving conditional expectation from $B^{\prime}\cap A_k$ onto $A^{\prime}\cap A_k$ is given by the following,
\begin{center}
$E^{B^{\prime}\cap A_k}_{A^{\prime}\cap A_k}(x)=\frac{1}{[A:B]_0}\sum_{i}\lambda_ix\lambda^*_i\,,\quad x \in B'\cap A_k.$
\end{center}
\end{lmma}

Recall that if $\{\lambda_i:1\leq i\leq n\}\subset A$ is a quasi-basis for $E_0$, then we have (see \cite{Watataniindex})
\begin{equation}\label{vip}
\sum_i \lambda_ie_1\lambda^*_i=1.
\end{equation}

\begin{crlre}[\cite{BakshiVedlattice}]\label{E-e}
Following the notation as in Lemma $\ref{f2}$, we have $E^{B^{\prime}\cap A_1}_{A^{\prime}\cap A_1}(e_1)=\delta^{-2}.$
\end{crlre}

\begin{rmrk}\rm
The reader should note that in \cite{Liunoncommutative} the authors have taken unnormalized trace $Tr_k$ on $N^{\prime}\cap M_k$, and thus $Tr_2(e_1)=1$. More precisely, $\tr_k(x)= {\delta}^{-(k+1)} Tr_k(x)$ for any $x\in B^{\prime}\cap A_k$.
\end{rmrk}

\subsection{Some useful inequalities}

For the convenience of the reader, we recall a few important inequalities as mentioned in \cite{Liunoncommutative}.
\begin{dfn}
For $x\in B^\prime\cap A_1$, we define the $p$-norm of $x$ for $1\leq p<\infty$ as follows~:
\begin{center}
$\lVert x\rVert_p=(\tr(|x|^p))^{\frac{1}{p}};$
\end{center}
and for $p=\infty$
\begin{center}
$\lVert x\rVert_{\infty}=\lVert x\rVert,$
\end{center}
where $\lVert.\rVert$ denotes the operator norm and $\tr$ denotes the Markov type trace on $B^\prime \cap A_1$ as in $\Cref{m1}$.
\end{dfn}

\begin{ppsn}[H\"{o}lder's Inequality]{\em \cite{Xunotes}}\label{holder} For any $x,y,z$ in $B^\prime \cap A_1$, we have the following,
\begin{enumerate}[$(i)$]
\item $|\tr(xy)|\leq ||x||_p||y||_q\,,\mbox{ where } 1\leq p\leq\infty\,,\,\,\frac{1}{p}+\frac{1}{q}=1;$
\item $|\tr(xyz)|\leq \|x\|_p\|y\|_q\|z\|_r\,,\mbox{ where } 1\leq p,q\leq\infty\,,\,\,\frac{1}{p}+\frac{1}{q}+\frac{1}{r}=1;$
\item $\|xy\|_r\leq \|x\|_p\|y\|_q\,,\mbox{ where } 1\leq p,q,r\leq\infty\,,\,\,\frac{1}{r}=\frac{1}{p}+\frac{1}{q}.$
\end{enumerate}
\end{ppsn}

\begin{ppsn}[\cite{Xunotes}]\label{norm}
For any $x$ in $B^\prime \cap A_1$ and $1\leq p<\infty$, we have
\begin{center}
$\|x\|_p=\sup\{|\tr(xy)|:y\in B^\prime\cap A_1,\|y\|_q\leq 1\},$
\end{center}
where $\frac{1}{p}+\frac{1}{q}=1$.
\end{ppsn}

\begin{ppsn}[\cite{kosaki}]\label{Interpolation}
Let $\mathcal{M}$ be a finite von Neumann algebra with a faithful normal tracial state $\tau$. Suppose that $T:\mathcal{M}\to\mathcal{M}$ is a linear map. If
\begin{center}
$\|Tx\|_{p_1}\leq K_1\|x\|_{q_1}\quad \mbox{ and }\quad\|Tx\|_{p_2}\leq K_2\|x\|_{q_2},$
\end{center}
then for any $\,\theta\in[0,1]$,
\begin{center}
$\|Tx\|_{p_\theta}\leq K_1^{1-\theta}K_2^{\theta}\,\|x\|_{q_\theta}$
\end{center}
where $\,\frac{1}{p_\theta}=\frac{1-\theta}{p_1}+\frac{\theta}{p_2}$ and $\frac{1}{q_\theta}=\frac{1-\theta}{q_1}+\frac{\theta}{q_2}$.
\end{ppsn}

%%%%%%%%%%%%%%%%%%%%%%%%%%%%%%%
%%%%%%%            Section 3            %%%%%%%%%%%%%
%%%%%%%%%%%%%%%%%%%%%%%%%%%%%%%

\newsection{Revisit of non-commutative Fourier theory}\label{Sec 3}

Analogous to subfactor theory, in \cite{BakshiVedlattice} the authors have provided a Fourier theory using Watatani’s notions of index and $C^*$-basic construction of certain inclusions of $C^*$-algebras. In this section we further investigate this Fourier theory and its properties which we shall use in the sequel. Throughout this section let $B\subset A$ denote an inclusion of simple $C^*$-algebras with a conditional expectation of finite Watatani index.

\subsection{The Rotation maps}\label{ro}

The Fourier transform of paragroups for a finite depth subfactor was first introduced by Ocneanu and as already mentioned in the introduction, it plays a major role in the development of the subfactor theory. More generally, for any extremal subfactor an explicit formula for the Fourier transform on the higher relative commutants was given by Bisch in (Def. $2.16$, \cite{Bi94}) (see also \cite{BiJo2} for many interesting results involving the Fourier transforms). The subtle difference between the Fourier theory for $C^*$-inclusion as in \cite{BakshiVedlattice} and that of subfactor theory lies in the fact that, unlike for finite factors, we neither have a tracial state on the $C^*$-algebra to begin with nor the ‘modular conjugation operator’.

\begin{dfn}
For each $k\geq 0$, the Fourier transform $\mathcal{F}_k:B^{\prime}\cap A_k\longrightarrow A^{\prime}\cap A_{k+1}$ is defined by the following,
\begin{center}
$\mathcal{F}_k(x)=\delta^{k+2}\,E^{B^{\prime}\cap A_{k+1}}_{A^{\prime}\cap A_{k+1}}(xe_{k+1}e_k\cdots e_2 e_1)\,.$
\end{center}
The inverse Fourier transform ${\mathcal{F}}^{-1}_k: A^{\prime}\cap A_{k+1}\longrightarrow B^{\prime}\cap A_k$ is defined by the following,
\begin{center}
$\mathcal{F}^{-1}_k(y)=\delta^{k+2}\,E_{k+1}(ye_1e_2\cdots e_ke_{k+1})\,.$
\end{center}
\end{dfn}

The meaning of ``inverse'' in the preceding definition is justified by the fact that $\mathcal{F}_k\circ\mathcal{F}^{-1}_k=\mathrm{id}_{A^\prime\cap A_{k+1}}$ and $\mathcal{F}^{-1}_k\circ\mathcal{F}_k=\mathrm{id}_{B^\prime\cap A_k}$ for all $k\geq 0$ (Proposition $3.2$ in \cite{BakshiVedlattice}). In this paper, we will be mainly interested in the case of $k=1$ and hence for simplicity, we denote the Fourier transform $\mathcal{F}_1$ by $\mathcal{F}$. Recall that both $\mathcal{F} \mbox{ and } \mathcal{F}^{-1}$ are isometries with respect to the norm given by $\|x\|_2 = (\tr(x^*x))^\frac{1}{2}$ (Theorem $3.5$ in \cite{BakshiVedlattice}).

We now revisit the rotation map defined in (Definition $3.7$ in \cite{BakshiVedlattice}) and derive a few more properties of it. Recall  the rotation map $\rho_{+}: B^\prime \cap A_1 \longrightarrow B^\prime \cap A_1$ defined by,
\begin{center}
$\rho_{+}(x)=(\mathcal{F}^{-1}((\mathcal{F}(x))^*))^*\,.$
\end{center}
It is known that $\rho_{+}$ is a unital involutive $*$-preserving anti-automorphism (Remark $3.11$ and Theorem $3.16$ in \cite{BakshiVedlattice}), and hence $\lVert x\rVert_\infty=\lVert\rho_{+}(x)\rVert_\infty$. It is also shown in \cite{BakshiVedlattice} that if the inclusion $B \subset A$ is irreducible, then $\rho_{+}$ is $\tr$-preserving.
  	
Analogous to $\rho_{+}$, we can define a  rotation operator $\rho_{-}: A^\prime \cap A_2 \longrightarrow A^\prime \cap A_2$ by the following,
\begin{center}
$\rho_{-}(w)=(\mathcal{F}((\mathcal{F}^{-1}(w))^*))^*.$
\end{center}

\begin{lmma}\label{rotationmaplemma}
We have $\,\rho_{-}=\mathcal{F}\circ \rho_{+}\circ \mathcal{F}^{-1}$. In other words, the diagram in \Cref{commutativediagram} commutes.
\begin{figure}[!h]
\begin{center}\begin{tikzcd}
B^\prime \cap A_1\arrow{r}{\mathcal{F}} \arrow[swap]{d}{\rho_{+}} & A^\prime \cap A_2 \arrow{d}{\rho_{-}} \\%
B^\prime \cap A_1\arrow{r}{\mathcal{F}}& A^\prime \cap A_2
\end{tikzcd}
\end{center}
\caption{Relation between $\rho_{+}$ and $\rho_{-}$}
\label{commutativediagram}
\end{figure}
\end{lmma}
\begin{prf}
Since $\rho_{+}$ is $*$-preserving, we have the following,
\begin{eqnarray}{\label{proof1}}
\mathcal{F}(x)^*=\mathcal{F}\circ \rho_{+}(x^*).
\end{eqnarray}
Now again from the definition,
\begin{eqnarray}\label{proof2}
 \rho_{-}(\mathcal{F}(x))=(\mathcal{F}(x^*))^*=\mathcal{F}\circ \rho_{+}(x),
 \end{eqnarray}
 where the last equality follows from Equation (\ref{proof1}). Thus, $\rho_{-}\circ \mathcal{F}=\mathcal{F}\circ \rho_{+}$.\qed
\end{prf}

Next we show that $\rho_{-}$ also satisfies properties  similar to $\rho_{+}$.

\begin{ppsn}\label{antihomo}
The map $\,\rho_{-}$ is a unital involutive $*$-preserving anti-automorphism.
\end{ppsn}
\begin{prf}
First we show that $\rho_{-}$ is $*$-preserving. To see this, for any $w\in A^{\prime}\cap A_2$ we observe the following,
\begin{eqnarray}\label{proof3}
w^*&=&(\mathcal{F}(\mathcal{F}^{-1}(w)))^*\nonumber\\
&=&\mathcal{F}\circ \rho_{+}((\mathcal{F}^{-1}(w))^*)
\end{eqnarray}
using \Cref{proof1}. Applying $\mathcal{F}^{-1}$ on both sides of the \Cref{proof3} we get the following,
\begin{eqnarray}\label{rotation1}
\mathcal{F}^{-1}(w^*))=\rho_{+}((\mathcal{F}^{-1}(w))^*).
\end{eqnarray}
Now apply $\rho_{+}$ on both sides of \Cref{rotation1} and use the fact that $\rho_{+}^2=\mbox{id}$ to get the following, 
\begin{eqnarray}\label{proof4}
(\mathcal{F}^{-1}(w))^*&=&\rho_{+}\circ \mathcal{F}^{-1}(w^*)\nonumber\\
&=&\mathcal{F}^{-1}\circ \rho_{-}(w^*)
\end{eqnarray}
using \Cref{rotationmaplemma}.
Finally, apply $\mathcal{F}$ on both sides of \Cref{proof4} and use the definition of $\rho_{-}$ to conclude that $\rho_{-}$ is $*$-preserving. The fact that $\rho_{-}^2 =\mbox{id}$ is an easy consequence of $\rho_{+}^2=\mbox{id}$ and Lemma \ref{rotationmaplemma}.

It remains to show that $\rho_{-}$ is a unital anti-homomorphism. If $\{\lambda_i:i\in I\}$ is a quasi-basis for $E_0$, then by \Cref{f2} we obtain the following,
\begin{eqnarray}\label{ma}
\rho_{-}(w) &=&(\mathcal{F}((\mathcal{F}^{-1}(w))^*))^*\nonumber\\
&=&\delta^3E^{B^\prime \cap A_2}_{A^\prime \cap A_2}(e_1e_2\mathcal{F}^{-1}(w))\nonumber\\
&=& \delta \sum_i \lambda_ie_1e_2\mathcal{F}^{-1}(w)\lambda_i^*\nonumber\\
&=& \delta^4 \sum_i \lambda_ie_1e_2E_2(we_1e_2)\lambda_i^* \,.
\end{eqnarray}
Since, $e_1e_2e_1={\delta}^{-2} e_1$ and $E_2(e_2)={\delta}^{-2}$, we have
\begin{center}
$\rho_{-}(1)= \delta^2 \sum_i \lambda_ie_1E_2(e_2)\lambda_i^*= \sum_i \lambda_ie_1\lambda_i^*=1.$
\end{center}
This last equality follows from Equation (\ref{vip}). As $\rho_{-}$ is $*$-preserving, to show that it is an anti-homomorphism, it is enough to show that, $\rho_{-}(w_1)\rho_{-}(w_2)^*=\rho_{-}(w_2^*w_1)$, for any $w_1,w_2\in A^{\prime}\cap A_2$. By Equation (\ref{ma}) and \Cref{pushdown} we finally get the following,
\begin{eqnarray}
\rho_{-}(w_1)\rho_{-}(w_2)^*&=& \delta^8 \sum_{i,j}\lambda_ie_1e_2 E_2(w_1e_1e_2)\lambda_i^* \lambda_jE_2(e_2e_1w_2^*)e_2e_1\lambda_j^* \nonumber\\
&=& \delta^6 \sum_{i,j}\lambda_ie_1E_1\circ E_2(w_1e_1\lambda_i^*\lambda_je_2e_1w_2^*)e_2e_1 \lambda_j^*\nonumber\\
&=& \delta^6 \sum_{i,j}\lambda_ie_1E_1\circ E_2(e_1\lambda_i^*\lambda_je_2e_1w_2^*w_1)e_2e_1\lambda_j^*~~~~~ \text{(Lemma 3.11 in \cite{KajiwaraWatatani})}\nonumber\\
&=& \delta^6 \sum_{i,j}\lambda_ie_1E_1\big(e_1E_2(\lambda_i^*\lambda_je_2e_1w_2^*w_1)\big)e_2e_1\lambda_j^*\nonumber\\
&=& \delta^4 \sum_{i,j}\lambda_ie_1E_2(\lambda_i^*\lambda_je_2e_1w_2^*w_1)e_2e_1\lambda_j^*~~~~\text{(by \Cref{pushdown})}\nonumber\\
&=& \delta^4 \sum_{i,j}\lambda_ie_1\lambda_i^*\lambda_jE_2(e_2e_1w_2^*w_1)e_2e_1 \lambda_j^*\nonumber\\
&=& \delta^4 \sum_j \lambda_jE_2(e_2e_1w_2^*w_1)e_2e_1 \lambda_j^*\nonumber\\
&=& \rho_{-}(w_1^*w_2)^*\nonumber\\
&=& \rho_{-}(w_2^* w_1)\nonumber
\end{eqnarray}
and this completes the proof.\qed
\end{prf}

\begin{ppsn}\label{tracepreserving}
If $B \subset A$ is irreducible, then $\rho_{-}$ is a $\tr$-preserving map on $A^\prime\cap A_2$, where $\tr$ on $A^\prime \cap A_2$ is the restriction of $\tr$ on $B^\prime\cap A_2$. 
\end{ppsn}
\begin{prf}
For $w\in A^\prime\cap A_2$, there exists a unique $x\in B^\prime\cap A_1$ such that $w=\mathcal{F}(x)$. By \Cref{rotationmaplemma}, we have
\begin{center}
$\tr(\rho_{-}(w))=\tr(\rho_{-}\circ\mathcal{F}(x))=\tr(\mathcal{F}\circ \rho_{+}(x))= \delta^3\,\tr\big(E^{B^{\prime}\cap A_2}_{A^{\prime}\cap A_2}(\rho_{+}(x)e_2e_1)\big)\,.$
\end{center}
Since $E^{B^{\prime}\cap A_2}_{A^{\prime}\cap A_2}$ is $\tr$-preserving, we get by \Cref{m1},
\begin{center}
$\tr(\rho_{-}(w))=\delta^3\,\tr(\rho_{+}(x)e_2e_1)=\delta\,\tr(e_1\rho_{+}(x))=\delta\,\tr(\rho_{+}(xe_1))\,.$
\end{center}
Here, the last equality follows from the fact that $\rho_{+}(e_1)=e_1$ and $\rho_{+}$ is an anti-homomorphism. Since $\rho_{+}$ is $\tr$-preserving, we immediately obtain the following,
\begin{center}
$\tr(\rho_{+}(xe_1))=\tr((xe_1))=\delta^2\,\tr(xe_2e_1)=\delta^2\,\tr(E^{B^{\prime}\cap A_2}_{A^{\prime}\cap A_2}(xe_2e_1))=\delta^{-1}\tr(\mathcal{F}(x))=\delta^{-1}\tr(w)\,.$
\end{center}
Therefore, we have $\tr(\rho_{-}(w))=\tr(w)$ as desired.\qed	
\end{prf}

\begin{crlre}\label{pnormcomparison1}
Let $1\leq p \leq \infty$. For an irreducible inclusion $B \subset A$, we have
\begin{center}
$\|x\|_p=\|\rho_{+}(x)\|_p\quad\quad\text{and}\quad \quad\|w\|_p =\|\rho_{-}(w)\|_p\,$
\end{center}
for $x\in B^\prime\cap A_1$ and $w\in A^{\prime}\cap A_2$.
\end{crlre}

\begin{rmrk}\label{pnormcomparison2}\rm
Note that \Cref{pnormcomparison1} need not be true for non-irreducible simple $C^*$-inclusions for $p\neq\infty$. For extremal $II_1$ factors, not necessarily irreducible, an easy pictorial calculation shows that $\rho_{+}$ (resp. $\rho_{-}$) being $\tr$-preserving is equivalent to \Cref{sphericality}, which in turn is equivalent to the sphericality, and hence \Cref{pnormcomparison1} holds.
	\begin{figure}[!h]
		\begin{center}
			\begin{tikzpicture}[scale=.55]
			\begin{scope}
			\begin{scope}
			\begin{scope}[shift={(0,1)}]
			\clip(-1,0)rectangle(1,1);
			\draw(0,0)circle(1);
			\end{scope}
			\begin{scope}[shift={(0,-1)}]
			\clip(-1,0)rectangle(1,-1);
			\draw(0,0)circle(1);
			\end{scope}
			\draw(1,1)--(1,-1)(-1,1)--(-1,-1);
			\end{scope}
			\begin{scope}[xscale=1.7,yscale=1.4]
			\begin{scope}[shift={(0,1)}]
			\clip(-1,0)rectangle(1,1);
			\draw(0,0)circle(1);
			\end{scope}
			\begin{scope}[shift={(0,-1)}]
			\clip(-1,0)rectangle(1,-1);
			\draw(0,0)circle(1);
			\end{scope}
			\draw(1,1)--(1,-1)(-1,1)--(-1,-1);
			\end{scope}
			\draw[fill=white](-4,-1)rectangle(0,1);
			\draw(-2,0)node{$x$};
			\draw(-3,1.5)node{*};
			\draw(2.5,0)node{$=$};
			\end{scope}
			\begin{scope}[shift={(5,0)},rotate=180]
			\begin{scope}
			\begin{scope}[shift={(0,1)}]
			\clip(-1,0)rectangle(1,1);
			\draw(0,0)circle(1);
			\end{scope}
			\begin{scope}[shift={(0,-1)}]
			\clip(-1,0)rectangle(1,-1);
			\draw(0,0)circle(1);
			\draw(0,-0.5)node{*};
			\end{scope}
			\draw(1,1)--(1,-1)(-1,1)--(-1,-1);
			\end{scope}
			\begin{scope}[xscale=1.7,yscale=1.4]
			\begin{scope}[shift={(0,1)}]
			\clip(-1,0)rectangle(1,1);
			\draw(0,0)circle(1);
			\end{scope}
			\begin{scope}[shift={(0,-1)}]
			\clip(-1,0)rectangle(1,-1);
			\draw(0,0)circle(1);
			\end{scope}
			\draw(1,1)--(1,-1)(-1,1)--(-1,-1);
			\end{scope}
			\draw[fill=white](-4,-1)rectangle(0,1);
			\draw(-2,0)node{$x$};
			\end{scope}
			\end{tikzpicture}
		\end{center}
		\caption{$\tr(x)=\tr\circ\rho_{+}(x)$}
		\label{sphericality}
	\end{figure}
\end{rmrk}

\subsection{Convolution}\label{co}

Using the Fourier transform, we can introduce a new multiplication structure on the relative commutant $B^\prime\cap A_1$ (resp. $A^{\prime}\cap A_2$), which we call the {\em convolution} product. This is defined formally below.

\begin{dfn}\label{coproduct}{\em \cite{BakshiVedlattice}}
The convolution product of two elements $x$ and $y$ in $B^{\prime}\cap A_1$, denoted by $x*y$, is defined as
\begin{center}
$x* y={\mathcal{F}}^{-1}\big(\mathcal{F}(y)\mathcal{F}(x)\big).$
\end{center}
Similarly, for any two elements $w, z\in A^{\prime}\cap A_2$ we define
\begin{center}
$w*z=\mathcal{F}\big({\mathcal{F}}^{-1}(z){\mathcal{F}}^{-1}(w)\big).$
\end{center}
\end{dfn}

Recall that (Lemma $3.20$, \cite{BakshiVedlattice}) the convolution $`*$' is associative. We now prove that it is well behaved with the adjoint operation.

\begin{ppsn}{\label{coproductadjoint}}
For $\,x,y\in B^\prime \cap A_1$, we have $\,(x* y)^*=(x^*)* (y^*)$. Similarly, for  $w, z\in A^{\prime}\cap A_2$, we have $\,(w* z)^*=(w^*)* (z^*)$.
\end{ppsn}
\begin{prf}
Let $\,x,y\in B^\prime \cap A_1$. Using \Cref{antihomo}, Equations (\ref{proof4}\,,\,\ref{proof1}) and \Cref{rotationmaplemma} we observe the following,
\begin{eqnarray}
(x* y)^* &=&(\mathcal{F}^{-1}(\mathcal{F}(y)\mathcal{F}(x)))^*\nonumber\\
&=&\mathcal{F}^{-1}(\rho_{-}((\mathcal{F}(y)\mathcal{F}(x))^*))\nonumber\\
&=& \mathcal{F}^{-1}(\rho_{-}((\mathcal{F}(x))^*(\mathcal{F}(y))^*))\nonumber\\
&=& \mathcal{F}^{-1}(\rho_{-}((\mathcal{F}(y))^*)\rho_{-}((\mathcal{F}(x))^*))\nonumber\\
&=& \mathcal{F}^{-1}(\rho_{-}\circ \mathcal{F}\circ \rho_{+}(y^*)\rho_{-}\circ \mathcal{F}\circ \rho_{+}(x^*))\nonumber\\
&=& \mathcal{F}^{-1}(\mathcal{F}(y^*)\mathcal{F}(x^*))\nonumber\\
&=& x^**y^*,\nonumber
\end{eqnarray}
which proves the first assertion, and the second assertion follows similarly.\qed
\end{prf}

We finally show that $\rho_{+}$ and $\rho_{-}$ are anti-multiplicative with respect to the convolution.

\begin{ppsn}\label{antimultiplicative}
For $x,y\in B^\prime\cap A_1$, we have $\rho_{+}(x * y)=\rho_{+}(y)*\rho_{+}(x)$. Similarly, for $w,z\in A^{\prime}\cap A_2$ one has $\,\rho_{-}(w * z)=\rho_{-}(z)* \rho_{+}(w)$.
\end{ppsn}
\begin{prf}
Observe that by Equations (\ref{proof2}\,,\,\ref{rotation1}) and \Cref{coproductadjoint}, we get
\begin{eqnarray}
\rho_{+}(y)* \rho_{+}(x) &=& \mathcal{F}^{-1}(\mathcal{F}\circ \rho_{+}(x)\mathcal{F}\circ\rho_{+}(y))\nonumber\\
&=& \mathcal{F}^{-1}((\mathcal{F}(x^*))^*(\mathcal{F}(y^*))^*)\nonumber\\
&=&\mathcal{F}^{-1}((\mathcal{F}(y^*)\mathcal{F}(x^*))^*)\nonumber\\
&=& \rho_{+}((\mathcal{F}^{-1}((\mathcal{F}(y^*)\mathcal{F}(x^*))^*)\nonumber\\
&=& \rho_{+}((x^* * y^*)^*)\nonumber\\
&=&\rho_{+}(x * y)\,,\nonumber
\end{eqnarray}
proving that $\rho_+$ is anti-multiplicative. A similar computation proves the result for $\rho_-$.\qed
\end{prf}

%%%%%%%%%%%%%%%%%%%%%%%%%%%%%%%
%%%%%%%%%            Section 4            %%%%%%%%%%%
%%%%%%%%%%%%%%%%%%%%%%%%%%%%%%%

\newsection{Fourier-theoretic inequalities}\label{Sec 4}

In this section we prove the Hausdorff-Young inequality and Young's inequality for inclusion $B\subset A$ of simple $C^*$-algebras which is not necessarily irreducible. We also provide various uncertainty principles on the second relative commutant of such an inclusion. In the non-commutative world, these Fourier-theoretic inequalities were first established in \cite{Liunoncommutative} for a finite index irreducible subfactor. We prove the non-commutative version of these inequalities for inclusion of simple $C^*$-algebras with a conditional expectation of index-finite type by finding the correct constants.
\smallskip

\textbf{Notation:} To avoid notational difficulty, for $x \in B^\prime \cap A_1$ and $w \in A^\prime \cap A_2$, we denote $\rho_{+}(x)$ and $\rho_{-}(w)$ by $\overline{x}$ and $\overline{w}$ respectively. In the case of inclusion of extremal $II_1$ factors, $\rho_{+}$ coincides with the 2-click rotation in the anti-clockwise direction, and hence it is consistent with the notations in \cite{Liunoncommutative}.
\smallskip

Define \[{\kappa}^{+}_0= \text{min}\big\{ \tr(p): p\in \mathcal{P}(B^{\prime}\cap A)\big\},\] 
\[{\kappa}^{-}_0= \text{min}\big\{ \tr(q): q\in \mathcal{P}(A^{\prime}\cap A_1)\big\}\] and 
\[{\kappa}_0=\sqrt{{\kappa}^{+}_0 {\kappa}^{-}_0}.\]
Recall that, $\delta=\sqrt{{[A:B]}_0}.$

\subsection{Hausdorff-Young inequality}\label{hau}

Goal of this subsection is to prove the following non-commutative analogue of the classical Hausdorff-Young inequality for inclusion of simple $C^*$-algebras $B\subset A$ with a conditional expectation of index-finite type.

\begin{thm}[Hausdorff-Young inequality]\label{HY}
Let $B\subset A$ be an inclusion of simple $C^*$-algebras with a conditional expectation of index-finite type. For any $x \in B^\prime \cap A_1$,
\begin{center}
$\Vert x \rVert_q \leq \lVert \mathcal{F}(x)\rVert_p\leq \Big(\dfrac{\delta}{\kappa_0}\Big)^{1-\frac{2}{p}}~\lVert x\rVert_q$
\end{center}
where, $2 \leq p \leq \infty$ and $\frac{1}{p}+\frac{1}{q}=1$.
\end{thm}

The main point is \Cref{maintheorem} which is instrumental in proving \Cref{HY}. To begin with, we prove a few useful lemmas.

\begin{lmma}\label{constant}
For $x\in B^{\prime}\cap A_1$, we have $\mathcal{F}(x){\mathcal{F}(x)}^*={\delta}^2\,E^{B^{\prime}\cap A_2}_{A^{\prime}\cap A_2}(xe_2x^*)$.
\end{lmma}
\begin{prf}
Let $\{\lambda_1,\cdots,\lambda_n\}\subset A$ be a \textit{quasi-basis} for the minimal conditional expectation $E_0$. Using \Cref{f2}, we observe that for $x\in B^{\prime}\cap A_1$ the following holds~:
\begin{eqnarray}
\mathcal{F}(x){\mathcal{F}(x)}^* &=& \delta^6\,E^{B^\prime \cap A_2}_{A^\prime \cap A_2}(xe_2e_1)E^{B^\prime \cap A_2}_{A^\prime \cap A_2}(e_1e_2x^*)\nonumber\\
&=& \delta^6\,E^{B^\prime \cap A_2}_{A^\prime \cap A_2}\left(xe_2e_1E^{B^\prime \cap A_2}_{A^\prime \cap A_2}(e_1e_2x^*)\right)\nonumber\\
&=& \delta^4\,E^{B^\prime \cap A_2}_{A^\prime \cap A_2}\big(\sum_i xe_2e_1 \lambda_i e_1 e_2 x^* \lambda_i^*\big)\nonumber\\
&=& \delta^4\,E^{B^\prime \cap A_2}_{A^\prime \cap A_2}\big(\sum_i x E_1(E_0(\lambda_i)e_1)e_2x^*\lambda_i^*\big)\nonumber\\
&=& \delta^2\,E^{B^\prime \cap A_2}_{A^\prime \cap A_2}\big(\sum_i xE_0(\lambda_i)e_2x^* \lambda_i^*\big)\nonumber\\
&=& \delta^2\,E^{B^\prime \cap A_2}_{A^\prime \cap A_2}\Big(xe_2x^*\big(\sum_iE_0(\lambda_i)\lambda_i^* \big)\Big)\nonumber\\
&=& \delta^2\,E^{B^\prime \cap A_2}_{A^\prime \cap A_2}(xe_2x^*)\nonumber
\end{eqnarray}
which completes the proof.\qed
\end{prf}

Recall  the following well-known result from the basic von Neumann algebra theory (see item 2.17 in \cite{stratila}, for instance).

\begin{lmma}\label{stratila}
If $\,0\leq a \leq 1$ and $p$ is a projection, then $0\leq a\leq p$ if and only if $a=ap$.
\end{lmma}

\begin{lmma}\label{normcomparison}
Suppose  that $\nu \in B^\prime \cap A_1$ is a non-zero partial isometry. Then, we have the following.
\begin{enumerate}[$(i)$]
\item $E_1(\nu^*\nu)\leq \dfrac{1}{{\kappa}_0^{+}} {\lVert \nu\rVert}_1.$
\item $E^{B^{\prime}\cap A_1}_{A^{\prime}\cap A_1}(\nu\nu^*)\leq \dfrac{1}{{\kappa}_0^{-}}{\lVert \nu\rVert}_1.$
\end{enumerate}
\end{lmma}
\begin{prf}
First note that $p=\nu^*\nu$ is a projection in  $B^\prime \cap A_1$ and so,
\begin{center}
$\lVert \nu \rVert_1=\tr|\nu|=\tr((\nu^*\nu)^{\frac{1}{2}})=\tr(\nu^*\nu).$
\end{center}
Take minimal projections $\{e_k\}_k \subset B^{\prime}\cap A$ such that $\sum_k e_k=1$. Thus, there are scalars $\alpha_k\geq 0$ such that $E_1(\nu^*\nu)=\sum_k \alpha_ke_k$. It follows that
\begin{center}
$tr(\nu^*\nu) =\sum_k\alpha_k\tr(e_k)\geq {\kappa}_0^{+}\sum_k\alpha_k\geq {\kappa}_0^{+}\sum_k\alpha_ke_k={\kappa}_0^{+} E_1(\nu^*\nu).$
\end{center}
This completes the proof for the first part. The proof for the other part is similar and we omit it.\qed
\end{prf}

\begin{lmma}\label{partialisometryclaim}
Suppose that $\nu \in B^\prime \cap A_1$ is a non-zero partial isometry. Then, we have
\begin{center}
$\nu e_2 \nu ^* \leq \dfrac{{\lVert \nu \rVert}_1}{{\kappa}_0^{+}} ~\nu \nu^*.$
\end{center}
\end{lmma}
\begin{prf}
Let $p=\nu \nu^*$  be the range projection of $\nu$. Now, using \Cref{normcomparison} we get
\begin{center}
$\nu e_2\nu^*.\nu e_2\nu^*=\nu E_1(\nu^*\nu)e_2\nu^*\leq \frac{{\lVert \nu\rVert}_1}{{\kappa}_0^{+}} \nu e_2\nu^*.$
\end{center}
After taking $\infty$-norm on both sides of the above inequality we have $\Big \lVert \frac{{\kappa}_0^{+}}{{\lVert \nu\rVert}_1} \nu e_2\nu^* \Big\rVert_{\infty}\leq 1$ and hence, $0 \leq \frac{{\kappa}_0^{+}}{{\lVert \nu\rVert}_1} \nu e_2\nu^* \leq 1$. Since $\nu$ is a partial isometry, we have $\big(\frac{{\kappa}_0^{+}}{\lVert\nu\rVert_1}\nu e_2 \nu ^*\big)\big(\nu\nu^*\big)=\frac{{\kappa}_0^{+}}{\lVert\nu\rVert_1}\nu e_2 \nu ^*$. Hence, the proof follows by Lemma \ref{stratila} with $a= \frac{{\kappa}_0^{+}}{\lVert\nu\rVert_1}\nu e_2 \nu ^*$ and $p=\nu \nu^*$.\qed
\end{prf}

\begin{ppsn}\label{maintheorem}
For $x \in B^\prime \cap A_1$, we have $$\,\lVert x\rVert _1\leq {\lVert \mathcal{F} (x)\rVert}_{\infty} \leq \dfrac{\delta}{\kappa_0} {\lVert x\rVert}_1.$$
\end{ppsn}
\begin{prf}
Note that by the Kadison-Schwarz inequality we have ${\lVert x\rVert }_1\leq {\lVert x\rVert}_2$, and we also know that ${\lVert x\rVert}_2\leq  {\lVert x\rVert}_{\infty}\,$. Since $\lVert \mathcal{F}(x)\rVert _2=\lVert x\rVert _2\,$, it follows that $\lVert x\rVert _1\leq {\lVert \mathcal{F} (x)\rVert}_{\infty}\,$. It remains to prove ${\lVert \mathcal{F} (x)\rVert}_{\infty} \leq \dfrac{\delta}{\kappa_0} {\lVert x\rVert}_1$. We first prove the inequality for a partial isometry and then appealing to rank-one decomposition as depicted in \cite{Liunoncommutative} proves the result for a general $x \in B^\prime \cap A_1$. Let $\nu$ be a partial isometry in $B^\prime \cap A_1$. By \Cref{constant} and \Cref{partialisometryclaim} we have the following,
\begin{eqnarray*}
\mathcal{F}(\nu){\mathcal{F}(\nu)}^*={\delta}^2\,E^{B^{\prime}\cap A_2}_{A^{\prime}\cap A_2}(\nu e_2\nu^*) \leq \dfrac{{\lVert \nu \rVert}_1}{{\kappa}_0^{+}} \delta^2 E^{B^{\prime}\cap A_1}_{A^{\prime}\cap A_1}(\nu\nu^*)\,.
\end{eqnarray*}
Now, using \Cref{normcomparison}(ii), we obtain the following,
\[\mathcal{F}(\nu){\mathcal{F}(\nu)}^*\leq \dfrac{\delta^2}{{\kappa}_0^2} {\lVert \nu\rVert}_1^2.\]
Applying $\|~.~\|_\infty$ on both sides and taking square root finishes the proof for $\nu$. For an arbitrary $x \in  B^\prime \cap A_1$, let $x=\sum_k \alpha_k \nu_k$ be the rank- one decomposition of $x$. Then, $\lVert x\rVert_1=\sum_k \alpha_k\lVert \nu_k \rVert_1$ and hence we have the following,
\[\|\f(x)\|_\infty\leq \sum_k \alpha_k\lVert \f(\nu_k) \rVert_\infty\leq \frac{\delta}{\kappa_0}\sum_k \alpha_k \|\nu_k\|_1 = \frac{\delta}{\kappa_0}\|x\|_1~,\]
and this completes the proof.\qed
\end{prf}
\smallskip

\textbf{Proof of \Cref{HY}:}
Using \Cref{maintheorem} and the fact that $\mathcal{F} \mbox{ and } \mathcal{F}^{-1}$ both are isometries with respect to the norm given by $\|x\|_2 = (\tr(x^*x))^\frac{1}{2}$, we have the following,
\[
{\lVert \mathcal{F} (x)\rVert}_{\infty}\leq \dfrac{\delta}{{\kappa}_0}{\lVert x\rVert}_1\quad\mbox{ and }\quad{\lVert \mathcal{F}(x)\rVert}_{2}={\lVert x\rVert}_2\,.
\]
The proof of the inequality $\lVert \mathcal{F}(x)\rVert_p\leq \Bigg(\dfrac{\delta}{\kappa_0}\Bigg)^{1-\frac{2}{p}}~\lVert x\rVert_q$  is now clear from \Cref{Interpolation} with $p_1=\infty $, $q_1=1$, $p_2=2$, $q_2=2$, $K_1=\dfrac{\delta}{\kappa_0}$, $K_2=1$ and $\theta=\frac{2}{p}$. The proof of $\lVert x\rVert_q\leq  \lVert \mathcal{F}(x)\rVert_p$ is similar.\qed

As a corollary of \Cref{HY}, we obtain the following Hausdorff-Young inequality for a finite index subfactor not necessarily irreducible and in particular, in the extremal case we recover the result in (Theorem $7.3$, \cite{Liunoncommutative}).

\begin{crlre}\label{miniHY}
Let $N\subset M$ be a subfactor with finite Jones index. Then, for any $x \in N^\prime \cap M_1$,
\[
\Vert x \rVert_q \leq \lVert \mathcal{F}(x)\rVert_p\leq \bigg(\dfrac{\sqrt{[M:N]_0}}{\kappa_0}\bigg)^{1-\frac{2}{p}}~\lVert x\rVert_q,
\]
where, $2 \leq p \leq \infty$ and $\frac{1}{p}+\frac{1}{q}=1$.
\end{crlre}

\subsection{Non-commutative uncertainty principles}\label{ncu}

Motivated by \cite{Liunoncommutative}, we prove the Donoho-Stark uncertainty principle and Hirschman-Beckner uncertainty principle for inclusion of simple $C^*$-algebras $B \subset A$ with a conditional expectation of index-finite type. Our proofs are essentially applications of \Cref{hau} and marginal modification of the proofs in \cite{Liunoncommutative} with revised constants. 

Recall that for $x \in B^\prime \cap A_1$, the range projection  $x$ is the smallest projection $l(x) \in B^\prime \cap A_1$ such that $l(x)x=x$. Now if $x= \sum_j \lambda_j\nu_j$ is the rank-one decomposition of $x$, then it is easy to see the following,
\begin{eqnarray}\label{rangeprojection}
l(x)=\sum_j \nu_j\nu_j^*.
\end{eqnarray} 
For $x \in B^\prime \cap A_1$, we denote $\mathcal{S}(x)=\tr(l(x))$.

\begin{thm}[Donoho-Stark uncertainty principle]\label{DS}
Consider an inclusion of simple $C^*$-algebras $B\subset A$ with a conditional expectation of index-finite type. For any non zero $x \in B^\prime \cap A_1$, we have
\[
\mathcal{S}(x)\mathcal{S}(\f(x))\geq \dfrac{\kappa^2_0
}{{[A:B]}_0}\,.
\] 
In particular, if $N\subset M$ is a subfactor with finite Jones index, then for any non zero $x\in N^{\prime}\cap M_1$ we have
\[
\mathcal{S}(x)\mathcal{S}(\f(x))\geq \dfrac{\kappa_0^2}{{[M:N]}_0}\,.
\]
\end{thm}
\begin{prf}
The proof is inspired by the proof of Theorem 5.2 in \cite{Liunoncommutative}. Let $x \in B^\prime \cap A_1$ and $\f(x)=\sum_j\lambda_j\nu_ j$ be the rank one decomposition of $\f(x)$. It is easy to see that $\mathcal{S}(\f(x))=\sum_j \|\nu_j\|_1$. By \Cref{holder} and \Cref{maintheorem} we have the following,
\begin{eqnarray}\label{donohostarkproof1}
\sup_j \lambda_j =\|\f(x)\|_\infty &\leq& \dfrac{\delta}{\kappa_0} \|x\|_1 \nonumber\\
&=& \dfrac{\delta}{\kappa_0}  \|l(x)x\|_1 \nonumber\\
&=& \dfrac{\delta}{\kappa_0}  \|x\|_2\|l(x)\|_2 \nonumber\\
&=& \dfrac{\delta}{\kappa_0}  \|\f(x)\|_2\|l(x)\|_2 \nonumber\\
&=& \dfrac{\delta}{\kappa_0}  \|\f(x)\|_2 (\mathcal{S}(x))^{\frac{1}{2}}\,.
\end{eqnarray} 
It is clear from the rank one decomposition of $\f(x)$ that $\|\f(x)\|_2=(\sum_j\lambda_j ^2 \|\nu_j\|_1)^{\frac{1}{2}}$. Hence, \Cref{donohostarkproof1} becomes,
\begin{eqnarray}
\sup_j \lambda_j &\leq& \dfrac{\delta}{\kappa_0}  \, \big(\sum_j\lambda_j ^2 \|\nu_j\|_1\big)^\frac{1}{2}(\mathcal{S}(x))^{\frac{1}{2}}\nonumber\\
&\leq & \dfrac{\delta}{\kappa_0}  \, (\sup_j \lambda_j )\big(\sum_j \|\nu _j \|_1\big)^\frac{1}{2}(\mathcal{S}(x))^{\frac{1}{2}}\nonumber\\
&=& \dfrac{\delta}{\kappa_0}  \, (\sup_j \lambda_j )(\mathcal{S}(\f(x)))^\frac{1}{2}(\mathcal{S}(x))^{\frac{1}{2}}\nonumber
\end{eqnarray}
which completes the proof.\qed
\end{prf}

Consider the continuous function $\eta: [0,\infty) \longrightarrow\mathbb{R}$ defined by
\begin{eqnarray}\label{etafunction}
\eta(t) &=& \begin{cases}
-t\log t  & \mbox{ if } t>0,\cr
0 &  \mbox{ if } t=0.
\end{cases}
\end{eqnarray}

\begin{dfn}[von Neumann entropy]\label{vonNeumannentropy}
For $x \in B^\prime \cap A_1$, the von Neumann entropy of $|x|^2$ is defined by the following,
\begin{center}
$H(|x|^2)=\tr(\eta(|x|^2)).$
\end{center}
\end{dfn}

\begin{thm}[Hirschman-Beckner uncertainty principle]\label{HBuncertainty}
Let $B\subset A$ be an inclusion of simple $C^*$-algebras with a conditional expectation of index-finite type. For any $x \in B^\prime \cap A_1$,
\[
\frac{1}{2}\big(H(|\f(x)|^2)+H(|x|^2)\big)\geq -\|x\|_2^2\,\Big(\log\bigg(\dfrac{\delta}{\kappa_0}\bigg) +\log \|x\|_2^2\Big)\,.
\]
In particular, if $\|x\|_2=1$, then we have
\[
 \frac{1}{2}\big(H(|\f(x)|^2)+H(|x|^2)\big)\geq-\log \bigg(\dfrac{\delta}{\kappa_0}\bigg)\,.
 \]
\end{thm}
\begin{prf}
The proof is a consequence of \Cref{HY} and the standard argument as in Theorem $5.5$ in \cite{Liunoncommutative}. However, we sketch the proof for completeness. Let $ 0\neq x \in B^\prime \cap A_1$ so that $\f(x) \neq 0$. By \Cref{HY} we have the following,
\begin{eqnarray}\label{hirsproof1}
\|\mathcal{F}(x)\rVert_p\leq \Big(\dfrac{\delta}{\kappa_0}\Big)^{1-\frac{2}{p}}~\lVert x\rVert_q
\end{eqnarray}
where, $2 \leq p \leq \infty$ and $\frac{1}{p}+\frac{1}{q}=1$. Consider the following function,
\[
f(p)=\log\|\f(x)\|_p-\log \|x\|_q-\log \bigg(\dfrac{\delta}{\kappa_0}\bigg)^{1-\frac{2}{p}}\,.
\]
By Equation (\ref{hirsproof1}), we have $f(p)\leq 0$. Now, since $\f$ is an isometry with respect to $\|.\|_2$, we have $f(2)=0$ and hence $f^\prime(2)\leq 0$. Hence, we obtain the following,
\[
\frac{\text{d}}{\text{d}p}\Big|_{p=2}\big(\|\mathcal{F}(x)\|_p^p\big)=-\frac{1}{2}H\big(|\mathcal{F}(x)|^2\big)\quad\mbox{and}\quad\frac{\text{d}}{\text{d}p}\Big|_{p=2}\big(\log\|\mathcal{F}(x)\|_p\big)=-\frac{1}{4}\log\|\mathcal{F}(x)\|_2^2\,-\frac{H(|\mathcal{F}(x)|^2)}{4\|\mathcal{F}(x)\|_2^2}\,.
\]
Similarly,
\[
\frac{\text{d}}{\text{d}p}\Big|_{p=2}\big(\log\|x\|_q\big)=\frac{1}{4}\log\|x\|_2^2+\frac{H(|x|^2)}{4\|x\|_2^2} \qquad\text{and}\qquad \frac{\text{d}}{\text{d}p}\Big|_{p=2}\Big(\log\bigg(\dfrac{\delta}{\kappa_0}\bigg)^{1-\frac{2}{p}}\Big)=\frac{1}{2}\log \bigg(\dfrac{\delta}{\kappa_0}\bigg)\,.
\]
Now, the above equations together with the facts $f^\prime(2)\leq 0$ and $\|\f(x)\|_2=\|x\|_2$ implies the following,
\begin{eqnarray}\label{hirsproof2}
-\frac{1}{4}\log\|x\|_2^2-\frac{1}{4}\frac{H(|\mathcal{F}(x)|^2)}{\|\mathcal{F}(x)\|_2^2}-\frac{1}{4}\log\|x\|_2^2-\frac{1}{4}\frac{H(|x|^2)}{\|x\|_2^2}-\frac{1}{2}\log \bigg(\dfrac{\delta}{\kappa_0}\bigg) \leq 0\,.
\end{eqnarray}
A rearrangement of Equation (\ref{hirsproof2}) completes the proof.\qed
\end{prf}

\begin{crlre}
Let $N\subset M$ be a subfactor with finite Jones index. Then, for any $x \in N^\prime \cap M_1$ we have
\[
\frac{1}{2}\big(H(|\f(x)|^2)+H(|x|^2)\big)\geq -\|x\|_2^2\,\big(\log\bigg(\dfrac{{\sqrt{[M:N]}_0}}{\kappa_0}\bigg) +\log\|x\|^2_2\big)\,.
\]
In particular, if $\|x\|_2=1$, then we have
\[
\frac{1}{2}\big(H(|\f(x)|^2)+H(|x|^2)\big)\geq-\log\bigg(\dfrac{{\sqrt{[M:N]}_0}}{\kappa_0}\bigg)\,.
\]
\end{crlre}

\subsection{Young's inequality}\label{Y}

Goal of this subsection is to prove the Young's inequality.  Throughout this subsection we fix an (not necessarily irreducible) inclusion of simple $C^*$-algebras $B\subset A$ with a conditional expectation of finite Watatani index.
 
\begin{thm}[Young's Inequality]\label{Younginequalitytheorem} Suppose $B\subset A$ is an inclusion of simple $C^*$-algebras with a conditional expectation of index-finite type. Then, for any $x,y\in B^\prime \cap A_1$, we have
\[
\|x* y\|_r\leq\dfrac{\delta}{\kappa^{+}_0} {\bigg(\frac{\|y\|_1}{\|\overline{y}\|_1}\bigg)}^{\frac{1}{r}}\|x\|_p\|\overline{y}\|_q
\]
where, $1\leq p,q,r\leq \infty$ and $\frac{1}{p}+\frac{1}{q}=\frac{1}{r}+1$.
\end{thm}

As a corollary, we prove Young's inequality for a subfactor which is not necessarily extremal. Recall that, a subfactor is extremal if and only if $[M:N]={[M:N]}_0$ and furthermore, in the extremal case for any $y\in N^{\prime}\cap M_1$ we have ${\lVert y\rVert}_1={\lVert \overline{y} \rVert}_1$ (see \Cref{pnormcomparison2}). We would like to mention that  in the extremal case we recover the following Young's inequality for spherical planar algebras as in (Theorem $7.6$, \cite{Liunoncommutative}).

\begin{crlre}[\cite{Liunoncommutative}]\label{sanatbaba}
If $N\subset M$ is an extremal subfactor with $[M:N]<\infty$, then for any $x,y\in N^{\prime}\cap M_1$ we have
\[
\|x* y\|_r\leq \dfrac{\sqrt{[M:N]}}{\kappa_0^+}~  \|x\|_p\|{y}\|_q
\]
where, $1\leq p,q,r\leq \infty$ and $\frac{1}{p}+\frac{1}{q}=\frac{1}{r}+1$.
\end{crlre}

To prove \Cref{Younginequalitytheorem} we start with proving a few results involving the Fourier transform $\mathcal{F}$ which will be crucially used.

\begin{lmma}\label{younginequalitylemma1}
For $x\in B^\prime \cap A_1$, we have $\mathcal{F}(x)e_1\mathcal{F}(x)^*=xe_2x^*$.
\end{lmma}
\begin{prf}
Let $\{\lambda_i:i\in I\}$ be a quasi-basis for $E_0\,$. Using \Cref{f2} we observe the following,
\begin{eqnarray}
\mathcal{F}(x)e_1\mathcal{F}(x)^* &=& \delta^6\,E^{B^\prime \cap A_2}_{A^\prime \cap A_2}(xe_2e_1)e_1E^{B^\prime \cap A_2}_{A^\prime \cap A_2}(e_1e_2x^*)\nonumber\\
&=& \delta^2 \sum_{i,j}\lambda_i xe_2e_1\lambda_i^* e_1\lambda_je_1e_2x^*\lambda_j^*\nonumber\\
&=& \delta^2 \sum_{i,j}\lambda_i x E_0(\lambda_i^*)E_0(\lambda_j)e_2e_1e_2x^*\lambda_j^*\nonumber\\
&=& \sum_{j}\Big(\sum_i\lambda_i  E_0(\lambda_i^*)\Big)E_0(\lambda_j) xe_2x^*\lambda_j^*\nonumber\\
&=& \sum_j E_0(\lambda_j) xe_2x^*\lambda_j^*\nonumber\\
&=& xe_2x^*\Big(\sum_j E_0(\lambda_j)\lambda_j^*\Big)\nonumber\\
&=& xe_2x^*\nonumber
\end{eqnarray}
which finishes the proof.\qed
\end{prf}

\begin{lmma}\label{younginequalitylemma2}
For $\,x,y\in B^\prime\cap A_1$, we have
\[
E_2\big(\mathcal{F}(x)y\mathcal{F}(x)^*\big)=\frac{1}{\delta}\big(y* (xx^*)\big)\,.
\]
\end{lmma}
\begin{prf}
Let $\{\lambda_i:i\in I\}$ be a quasi-basis for $E_0$. Using \Cref{f2} we observe the following,
\begin{eqnarray}\label{lemmaproof}
E_2(\mathcal{F}(x)y\mathcal{F}(x)^*) &=& \delta^6\,E_2\big(E^{B^\prime \cap A_2}_{A^\prime \cap A_2}(xe_2e_1)yE^{B^\prime \cap A_2}_{A^\prime \cap A_2}(e_1e_2x^*)\big)\nonumber\\
&=& \delta^2 \sum_{i,j}E_2(\lambda_ixe_2e_1\lambda_i^*y\lambda_je_1e_2x^*\lambda_j^*) \nonumber\\
&=& \sum_{i,j}\lambda_ixE_1(e_1\lambda_i^*y\lambda_je_1)x^*\lambda_j^*\,.
\end{eqnarray}
Since $y \in B^\prime \cap A_1$, we can write $y=y_0e_1y_1$ for some $y_0,y_1 \in A$. Then, from Equation (\ref{lemmaproof}) and again using \Cref{f2}, we have the following,
\begin{eqnarray}
 E_2(\mathcal{F}(x)y\mathcal{F}(x)^*)&=&\sum_{i,j}\lambda_ixE_1(e_1\lambda_i^*y_0e_1y_1\lambda_je_1)x^*\lambda_j^*\nonumber\\
 &=& \sum_{i,j} \lambda_iE_0(\lambda_i^*y_0)xx^*E_1(e_1)E_0(y_1\lambda_j)\lambda_j^*\nonumber\\
 &=& \sum_{i}\lambda_ixx^*E_1(E_0(\lambda_i^*y_0)e_1y_1)\nonumber\\
 &=& \sum_i \lambda_ixx^*E_1(e_1\lambda_i^*y)\nonumber\\
 &=&\sum_{i,j} \lambda_ixx^*E_1(e_1\lambda_i^*\lambda_j E_0(\lambda_j^*)y)\nonumber\\
 &=& \delta^2  \sum_{i,j} \lambda_ixx^*E_1(e_1\lambda_i^*\lambda_jy)E_0(\lambda_j^*)E_2(e_2)\nonumber\\
 &=& \delta^4  \sum_{i,j}E_2(\lambda_ixx^*E_1(e_1\lambda_i^*\lambda_jy)e_2E_0(\lambda_j^*)e_1e_2)\nonumber\\
 &=&\delta^4  \sum_{i,j}E_2(\lambda_ixx^*e_2e_1\lambda_i^*\lambda_jye_2e_1\lambda_j^*e_1e_2)\nonumber\\
 &=& \delta^8\,E_2(E^{B^\prime \cap A_2}_{A^\prime \cap A_2}(xx^*e_2e_1)E^{B^\prime \cap A_2}_{A^\prime \cap A_2}(ye_2e_1)e_1e_2)\nonumber\\
 &=& \frac{1}{\delta}\,\mathcal{F}^{-1}(\mathcal{F}(xx^*)\mathcal{F}(y))\nonumber\\
 &=& \frac{1}{\delta}\,y * (xx^*)\nonumber
\end{eqnarray}
which completes the proof.\qed
\end{prf}

As a corollary, we prove the following Schur product theorem. In the planar algebraic language this was first noticed by Liu (Theorem $4.1$ in \cite{liu2016exchange}).

\begin{crlre}(Schur product theorem)\label{schurproduct}
If $x,y \in  B^\prime \cap A_1$ are positive, then $x* y$ is positive.
\end{crlre}
\begin{prf}
Let $x=aa^*,\,y=bb^*$ for some $a,b \in B^\prime \cap A_1$. Then, by \Cref{younginequalitylemma2} and in view of the fact that $E_2$ is a positive map it is now easy to see that $x*y$ is positive. Indeed,
\begin{eqnarray}
x*y=(aa^*)* (bb^*) &=& \delta\,E_2\big(\mathcal{F}(b)aa^* \mathcal{F}(b)\big)\nonumber\\
&=& \delta\,E_2\big(\mathcal{F}(b)a(\mathcal{F}(b)a)^*\big)\nonumber\\
&\geq& 0\,.\nonumber
\end{eqnarray}\qed
\end{prf}

\begin{rmrk}\rm
We would like to remark that as a consequence of \Cref{younginequalitylemma2}, we can prove that the `coproduct' on $B^\prime \cap A_1$ is well behaved with adjoints. However, at present we are not sure whether these two notions are equivalent. To see this, it is enough to take $x, y\in B^\prime \cap A_1$ such that $y\geq 0$. We write $y=bb^*$ for some $b \in B^\prime \cap A_1$. Then by \Cref{younginequalitylemma2} and since $E_2$ is $*$ - preserving, we have:
\begin{eqnarray}
(x * y)^*&=& (x*(bb^*))^*\nonumber\\
&=& \delta^{-1}E_2(\f(b)x\f(b)^*)^*\nonumber\\
&=& \delta^{-1}E_2(\f(b)x^*f(b)^*)\nonumber\\
&=& (x^* * (bb^*))\nonumber\\
&=& (x^**y^*).\nonumber
\end{eqnarray}
\end{rmrk}

Next we prove a Frobenius reciprocity type result as follows. 

\begin{crlre}\label{coproductapplication}
For $x,y,z \in B^\prime \cap A_1$, we have
\[
\tr((x* y)z)=\tr(x(z* \overline{y})
\]
\end{crlre}
\begin{prf}
First we claim that $\tr\circ E_2=\tr_2$, where $\tr_2$ is the Markov type trace on $B^\prime \cap A_2$ (see \Cref{m1}). Observe the following for $x \in B^\prime \cap A_2$~:
\begin{eqnarray}
\tr_1 \circ E_2(x) &=& E_0 \circ E_1|_{B^\prime \cap A_1}(E_2(x))\nonumber\\
&=& E_0 \circ E_1 \circ E_2|_{B^\prime \cap A_2}(x)\nonumber\\
&=& \tr_2(x)\nonumber
\end{eqnarray}
which finishes the proof of the claim.

Recall that by \Cref{proof1}, we have $(\mathcal{F}(x))^*=\mathcal{F}\big(\overline{x^*}\big)$. Now to prove the statement, assume that $y$ is positive so that $y=bb^*$ for some $b \in B^\prime \cap A_1$. Then, by \Cref{younginequalitylemma2} we have the following,
\begin{eqnarray}
\tr(x * y)z) &=& \tr((x* (bb^*))z)\nonumber\\
&=& \delta\,\tr(E_2(\mathcal{F}(b)x\mathcal{F}(b)^*z)\nonumber\\
&=& \delta\,\tr_2(\mathcal{F}(b)x\mathcal{F}(b)^*z)\nonumber\\
&=& \delta\,\tr_2(x \mathcal{F}(b)^*z\mathcal{F}(b))\nonumber\\
&=& \delta\,\tr(xE_2(\mathcal{F}(\overline{b^*})z(\mathcal{F}(\overline{b^*}))^*))\nonumber\\
&=& \tr(x(z *(\overline{b^*}(\overline{b^*})^*)))\nonumber\\
&=& \tr(x(z * \overline{y}))\nonumber
\end{eqnarray}
where, the last equation follows from \Cref{antihomo}.\qed
\end{prf}

The following lemma is crucial in proving the Young's inequality.

\begin{lmma}\label{young1}
For any $x,y\in B^\prime \cap A_1$, we have
\[
\|x* y\|_\infty\leq\dfrac{\delta}{\kappa^{+}_0} \,{\|x\|_\infty\|\overline{y}\|_1}\quad\mbox{and}\quad \|y* x\|_\infty\leq \dfrac{\delta}{\kappa^{+}_0} \,{\|x\|_\infty\|y\|_1}\,.
\]
\end{lmma}
\begin{prf}
First we prove that for $w\in A^{\prime}\cap A_2$ we have the following,
\begin{eqnarray}\label{younginequalitylemma3}
\mathcal{F}^{-1}(w){\mathcal{F}^{-1}(w)}^* &=& {\delta}^2\,E_2(we_1w^*)\,.
\end{eqnarray}
To see this, using \Cref{pushdown} we observe the following,
\begin{eqnarray}
\mathcal{F}^{-1}(w){\mathcal{F}^{-1}(w)}^*&=& \delta^6\, E_2(we_1e_2)E_2(e_2e_1w^*)\nonumber\\
&=& \delta^6\, E_2(we_1e_2E_2(e_2e_1w^*))\nonumber\\
&=& \delta^4\, E_2(we_1e_2e_1w^*)\nonumber\\
&=& \delta^2\,E_2(we_1w^*)\,.\nonumber
\end{eqnarray}
Next, suppose $\nu \in B^\prime \cap A_1$ is a partial isometry. For $x \in B^\prime \cap A_1$, using Equation (\ref{younginequalitylemma3}) we observe the following,
\begin{eqnarray}
(\nu * x)(\nu * x)^* &=& \mathcal{F}^{-1}(\mathcal{F}(x)\mathcal{F}(\nu))(\mathcal{F}^{-1}(\mathcal{F}(x)\mathcal{F}(\nu)))^* \nonumber\\
&=& \delta ^2\,E_2(\mathcal{F}(x)\mathcal{F}(\nu)e_1(\mathcal{F}(\nu))^*(\mathcal{F}(x))^*) \nonumber\\
&=& \delta ^2\,E_2(\mathcal{F}(x)\nu e_2 \nu^* \mathcal{F}(x)^*)\nonumber \,.
\end{eqnarray}
Now, by \Cref{partialisometryclaim} we have $\nu e_2 \nu ^* \leq \dfrac{\lVert\nu\rVert_1}{{\kappa}_0^{+}} \nu \nu^*$. Then, using \Cref{younginequalitylemma2}, we observe that
\begin{eqnarray}\label{younginequalityproof}
(\nu * x)(\nu * x)^* &\leq & \delta ^2 \frac{\lVert\nu\rVert_1}{{\kappa}_0^{+}} E_2(\mathcal{F}(x)\nu \nu^*\mathcal{F}(x)^*)\nonumber\\
&=& \delta \frac{\lVert\nu\rVert_1}{{\kappa}_0^{+}} (\nu\nu^*)*(xx^*)\nonumber\\
& \leq & \delta \frac{\lVert\nu\rVert_1}{{\kappa}_0^{+}} {\lVert x\rVert}^2_{\infty}(\nu\nu^*)*1.
\end{eqnarray}
On the other hand, we note that 
\begin{eqnarray}\label{youngirredproofone}
(\nu\nu^*)* 1 &=& \f^{-1}(\f(1)\f(\nu\nu^*)) \nonumber\\
&=& \delta ^3 E_2(\f(1)\f(\nu\nu ^*)e_1e_2) \nonumber\\
&=& \delta^4 E_2(e_2\f(\nu\nu ^*)e_1e_2) \nonumber\\
&=& \delta ^7 E_2(E^{B^{\prime}\cap A_2}_{A^{\prime}\cap A_2}(e_2\nu\nu ^* e_2e_1)e_1e_2) \nonumber\\
&=& \delta ^7 E_2(E^{B^{\prime}\cap A_2}_{A^{\prime}\cap A_2}(E_1(\nu\nu ^*)e_2e_1)e_1e_2) \,.
\end{eqnarray}
Now, using \Cref{E-e}, \Cref{normcomparison} and \Cref{youngirredproofone} we conclude that
$$(\nu\nu^*)* 1 \leq \frac{\delta}{{\kappa}_0^{+}}{\lVert \nu \rVert}_1.$$
Therefore, from \Cref{younginequalityproof} it follows that
$$ (\nu * x)(\nu * x)^* \leq \bigg(\dfrac{\delta}{{{\kappa}_0}^{+}}\bigg)^2{\lVert \nu\rVert}^2_1 {\lVert x\rVert}^2_{\infty}.$$
Thus we obtain, 
\begin{eqnarray}\label{younginequalityproof4}
\lVert \nu * x\rVert_\infty \leq \dfrac{\delta}{\kappa^{+}_0}\,\lVert\nu\rVert_1~ \lVert x\rVert_\infty\,.
\end{eqnarray}
Now, let $y \in B^\prime \cap A_1$ be arbitrary and $y=\sum_k \lambda_k \nu_k$ be the rank one decomposition of $y$. Then, by Equation (\ref{younginequalityproof4}) we have the following,
\begin{eqnarray}
\lVert y * x\rVert_\infty &\leq & \sum_k \lambda_k\,\lVert \nu * x\rVert_\infty\nonumber\\
&\leq &\dfrac{\delta}{\kappa^{+}_0} \,\sum_k \lambda_k\,\lVert\nu\rVert_1~ \lVert x\rVert_\infty\nonumber\\
&=& \dfrac{\delta}{\kappa^{+}_0} \,\lVert x\rVert_\infty\,\lVert y\rVert_1,\nonumber
\end{eqnarray}
and by \Cref{antimultiplicative} we have the following,
\begin{eqnarray}
\lVert x * y\rVert_\infty=\lVert \overline{x * y}\rVert_\infty&=&\lVert \overline{y}* \overline {x}\rVert_\infty\nonumber\\
&\leq & \dfrac{\delta}{\kappa^{+}_0} \,\lVert \overline{x}\rVert_\infty~\lVert \overline{y}\rVert_1\nonumber\\
&= & \dfrac{\delta}{\kappa^{+}_0} \,\lVert x\rVert_\infty~\lVert \overline{y}\rVert_1\,.\nonumber
\end{eqnarray}
Here $\|x\|_\infty=\|\overline{x}\|_\infty$, since the map $\rho_{+}$ is a unital anti-homomorphsim. This completes the proof.\qed
\end{prf}

\begin{lmma}\label{young2}
For any $x,y\in B^\prime \cap A_1$, we have
\[
\|x* y\|_1\leq \dfrac{\delta}{\kappa^{+}_0} {\|x\|_1\|y\|_1}\,.
\]
\end{lmma}
\begin{prf}
For any $x,y\in B^\prime \cap A_1$, using \Cref{norm}, \Cref{coproductapplication}, \Cref{holder} and \Cref{young1} respectively, we get the following,
\begin{eqnarray}
\|x* y\|_1=\sup_{\|z\|_\infty=1}|\tr((x* y)z)| &=&  \sup_{\|z\|_\infty=1}|\tr(x(z* \overline{y}))|\nonumber\\
&\leq&  \|x\|_1 \|z*\overline{y}\|_\infty\nonumber\\
&\leq& \dfrac{\delta}{\kappa^{+}_0} \|x\|_1 {\|\overline{\overline{y}}\|_1}\nonumber\\
&=& \dfrac{\delta}{\kappa^{+}_0}  \|x\|_1\|y\|_1\,.\nonumber
\end{eqnarray}\qed
\end{prf}

\begin{lmma}\label{young3}
For any $x,y\in B^\prime \cap A_1$, we have
\[
\|x* y\|_p\leq \dfrac{\delta}{\kappa^{+}_0}\bigg(\dfrac{\|y\|_1}{\|\overline{y}\|_1}\bigg)^{\frac{1}{p}} {\|x\|_p\|\overline{y}\|_1}~~\quad\mbox{and}\quad \|y* x\|_p\leq \dfrac{\delta}{\kappa^{+}_0}  {\|x\|_p\|y\|_1}
\]
where $1\leq p\leq \infty$.
\end{lmma}
\begin{prf}
For fixed $y \in B^\prime \cap A_1$, define $T_y: B^\prime \cap A_1 \longrightarrow B^\prime \cap A_1$ by
\begin{center}
$T_y(x)=x * y.$
\end{center}
Clearly $T_y$ is linear. Now, \Cref{young1} and \Cref{young2} respectively implies the following,
\begin{center}
$\|T_y(x)\|_\infty=\|x * y\|_\infty \leq \dfrac{\delta}{\kappa^{+}_0} \|x\|_\infty \|\overline{y}\|_1\,,$
\end{center}
and
\begin{center}
$\|T_y(x)\|_1=\|x * y\|_1 \leq \dfrac{\delta}{\kappa^{+}_0} \|x\|_1 \|y\|_1\,.$
\end{center}
Applying \Cref{Interpolation} with $p_1=\infty$, $p_2=1$, $q_1=\infty$, $q_2=1$, $\theta=\dfrac{1}{p}$, $K_1=\dfrac{\delta}{\kappa^{+}_0} \|\overline{y}\|_1$ and  $K_2= \dfrac{\delta}{\kappa^{+}_0} \|y\|_1$ we get  the first part. For the second part, for fixed $y$ we define $T_y(x)=y * x$. Then, a similar proof as above implies the result.\qed 
\end{prf}

\begin{lmma}\label{young4}
For any $x,y\in B^\prime \cap A_1$, we have
\[
\|x* y\|_\infty\leq \dfrac{\delta}{\kappa^{+}_0} \,{\|x\|_p\|\overline{y}\|_q}
\]
where $1\leq p \leq \infty$ and $\frac{1}{p}+\frac{1}{q}=1$.
\end{lmma}
\begin{prf}
Let $x * y=\sum_k \lambda_k \nu_k$ be the rank-one decomposition of $x * y$. Then,
\begin{eqnarray}\label{youngeq1}
\|x* y\|_\infty=\sup_k\,\lambda_k=\sup_k\,\frac{tr((x* y)v_k^*)}{tr(|v_k|)}.
\end{eqnarray}
Using \Cref{coproductapplication}, \Cref{holder} and \Cref{young3} respectively we see that
\begin{eqnarray}\label{youngeq2}
|\tr((x* y)v_k^*)| &=& |\tr(x(\nu_k^* * \overline{y}))\nonumber\\
&\leq & \|x\|_p\|\nu_k^* * \overline{y}\|_q \nonumber\\
&\leq & \dfrac{\delta}{\kappa^{+}_0} \,\|x\|_p \|\overline{y}\|_q \|\nu_k^*\|_1\nonumber\\
& \leq & \dfrac{\delta}{\kappa^{+}_0} \,\|x\|_p \|\overline{y}\|_q\,\tr(|\nu_k^*|)\,.
\end{eqnarray}
The proof is now clear from Equations (\ref{youngeq1} and  \ref{youngeq2}).\qed
\end{prf}
\smallskip

\textbf{Proof of \Cref{Younginequalitytheorem}~:} Fix $ x \in B^\prime \cap A_1$ and define $T_x: B^\prime \cap A_1 \to B^\prime \cap A_1$ by $T_x(y)=x *\overline{y}$. For $x,y\in B^\prime \cap A_1$, thanks to \Cref{young3} and \Cref{young4}, we have the following,
\begin{center}
$\|T_x(\overline{y})\|_p=\|x* y\|_p\leq \dfrac{\delta}{\kappa^{+}_0}\,{\|x\|_p\|\overline{y}\|_1}^{-\frac{1}{p}}\|y\|_1^{\frac{1}{p}}\|\overline{y}\|_1$
\end{center}
and
\begin{center}
$\|T_x(\overline{y})\|_\infty=\|x* y\|_\infty\leq \dfrac{\delta}{\kappa^{+}_0}\,{\|x\|_p\|\overline{y}\|_{\frac{1}{\frac{1}{q}-\frac{1}{r}}}}\,.$
\end{center}
The proof is now clear by the \Cref{Interpolation} with $p_1=p$, $p_2=\infty$, $q_1=1$, $q_2=\dfrac{1}{\frac{1}{q}-\frac{1}{r}}$, $K_1=\dfrac{\delta}{\kappa^{+}_0} \|x\|_p\|\overline{y}\|_1^{-\frac{1}{p}}\|y\|_1^{\frac{1}{p}}$, $K_2=\dfrac{\delta}{\kappa^{+}_0} \|x\|_p$ and $\theta=1-\dfrac{p}{r}\,$.\qed

\begin{rmrk}\rm
We remark that if $B\subset A$ is an irreducible inclusion of simple $C^*$-algebras with a conditional expectation of index-finite type, then the quadruple $(B^\prime \cap A_1,\tr,*,\rho_{+})$ forms a Frobenius $\delta$-algebra. We refer the reader to \cite{liuquantumfrobenius} for the definition of a Frobenius $k$-algebra. We also send the reader to \cite{liuquantumsmooth} for related notions.
\end{rmrk}

%%%%%%%%%%%%%%%%%%%%%%%%%%%%%%%%%%%%%%%%%%%%%%
%%%%%%%%%%%%%%%%%%    Irreducible    %%%%%%%%%%%%%%%%%%%%%
%%%%%%%%%%%%%%%%%%%%%%%%%%%%%%%%%%%%%%%%%%%%%%

\section{Appendix}

In this appendix, we discuss two examples to illustrate our results. We skip the proofs as they are routine verification. In order to investigate theory concerning Fourier transform for inclusion of simple $C^*$-algebras, these two examples can be considered as model examples to test new theories.

\subsection{Fourier transform for noncommutative torus}\label{NCT}

Let $\theta$ be an irrational number and consider the universal $C^*$-algebra $\mathscr{A}_\theta$, called the noncommutative torus, generated by two unitary elements $U$ and $V$ satisfying $UV=e^{-2\pi i\theta}VU$. It has a unital dense subalgebra $\mathbb{T}_\theta$ given by the following,
\begin{center}
$\mathbb{T}_\theta:=\big\{a=\displaystyle{\sum_{m,n\in\bbz}}a_{m,n}U^mV^n:\{a_{m,n}\}\in\mathcal{S}(\bbz^2)\big\},$
\end{center}
where $\mathcal{S}(\bbz^2)$ is the space of rapidly decreasing double sequences. The $C^*$-algebra $\mathscr{A}_\theta$ is equipped with a distinguished faithful {\em tracial state}, given on the dense subalgebra $\mathbb{T}_\theta$ by $\tau(a)=a_{0,0}$ and extends to $\mathscr{A}_\theta$ by continuity. We refer to (Chapter $6$, Section $3$ in \cite{connesgeometry}) for these facts. Moreover, $\mathscr{A}_\theta$ is a simple $C^*$-algebra since $\theta$ is irrational.

Let $k\geq 2$ be any natural number. Let us consider the unital $C^*$-subalgebra $\mathcal{B}_\theta$ of $\mathscr{A}_\theta$ generated by $U^k$ and $V$. By the universality and simplicity of $\mathscr{A}_{k\theta}$, it follows that $\mathcal{B}_\theta$ is canonically isomorphic to $\mathscr{A}_{k\theta}$. Assume further that $\theta$ is not an algebraic number of degree $2$. Then, the Watatani index $[\mathscr{A}_\theta:\mathcal{B}_\theta]_0$ is equal to $k$ (Page $112$ in \cite{Watataniindex}). Observe that $\mathcal{B}_\theta$ is nothing but the fixed point subalgebra of $\mathscr{A}_\theta$ under the $\bbz_k$ action given by $\overline{m}.U=e^{2\pi im/k}U$ and $\overline{m}.V=V$ for all $\overline{m}\in\bbz_k$. This says that the inclusion $\mathcal{B}_\theta\subset\mathscr{A}_\theta$ is in fact $\mathscr{A}_\theta^{\bbz_k}\subset\mathscr{A}_\theta$, and hence the basic construction is $\mathscr{A}_\theta\rtimes\bbz_k$. Since $\bbz_k$ is abelian, it then turns out that the Fourier transform $\mathcal{F}$ simply becomes the Fourier transform from the group algebra $\bbc\bbz_k$ onto $\bbc\widehat{\bbz}_k\cong\bbc\bbz_k$, and hence given by the $k\times k$ Fourier matrix.
\smallskip

However, in this example one can take the pedestrian way to find the explicit form of the Fourier transform and the rotation maps without invoking any result involving the crossed product. Observe that (see Proposition $(2.2.11,2.2.12)$ in \cite{Watataniindex}) we may use the GNS construction to realize the $C^*$-basic construction in this case. We only mention the intermediate steps and leave the detail to the interested reader for verification. For notational simplicity, we denote $\mathcal{B}\subset\mathcal{A}$ to mean the inclusion $\mathcal{B}_\theta\subset\mathscr{A}_\theta$. The GNS Hilbert space $L^2(\mathscr{A}_\theta,\tau)$ is isomorphic to $\ell^2(\bbz^2)$ via the identification $U^mV^n\longmapsto e_{m,n}$. Define $k$-many mutually orthogonal projections $p_r\in\mathscr{B}\left(\ell^2(\bbz^2)\right),\,0\leq r\leq k-1,$ by $p_r:e_{m,n}\mapsto e_{m,n}$ if $m\in k\bbz+r$ and $0$ otherwise. Considering the unital $C^*$-subalgebra $\mathcal{A}_1$ of $\mathscr{B}\left(\ell^2(\bbz^2)\right)$ generated by $U,V$ and $p_0$ (here $\tau(p_0)=1/k$) one gets the basic construction $\mathcal{B}\subset\mathcal{A}\subset\mathcal{A}_1$. Using the decomposition $\mathcal{A}_1=\bigoplus_{r=0}^{k-1}\,p_r\mathscr{A}_\theta$ as inner-product space, the GNS Hilbert space $L^2(\mathcal{A}_1,\tau)$ is isomorphic to $\bbc^k\otimes\ell^2(\bbz^2)$ by the following map,
\begin{center}$
p_0U^{m_0}V^{n_0}\,\oplus\ldots\ldots\oplus\,p_{k-1}U^{m_{k-1}}V^{n_{k-1}}\longmapsto\frac{1}{\sqrt{k}}(e_{m_0,n_0}\,,\ldots,\,e_{m_{k-1},n_{k-1}})\,.$
\end{center}
The following orthogonal projection $q:\bbc^k\otimes\ell^2(\bbz^2)\longrightarrow\bbc^k\otimes\ell^2(\bbz^2)$
\begin{center}$
q:(e_{m_0,n_0}\,,\ldots,\,e_{m_{k-1},n_{k-1}})\longmapsto \frac{1}{k}\big(\sum_{r=0}^{k-1}e_{m_r,n_r}\,,\ldots,\,\sum_{r=0}^{k-1}e_{m_r,n_r}\big)$
\end{center}
has range $\ell^2(\bbz^2)\cong L^2(\mathcal{A},\tau)$, and we obtain the basic construction tower of simple $C^*$-algebras $\mathcal{B}\,\subset\,\mathcal{A}\,\subset^{\,p}\,\mathcal{A}_1\,\subset^{\,q}\,\mathcal{A}_2$, where $\mathcal{A}_2=C^*\{\mathcal{A}_1,q\}\subseteq\mathscr{B}\left(\bbc^k\otimes\ell^2(\bbz^2)\right)=M_k(\bbc)\otimes\mathscr{B}\left(\ell^2(\bbz^2)\right)$. It follows that $\mathcal{A}_2=M_k(\bbc)\otimes\mathscr{A}_\theta$. For a subset $S\subseteq M_k(\bbc)$, we denote by $\mbox{Alg}\{S\}$ the subalgebra generated by $S$ and $S^*=\{x^*:x\in S\}$ in $M_k(\bbc)$. Let $C_k$ denote the permutation matrix $E_{1,k}+\sum_{i=1}^{k-1}E_{i+1,i}$ in $M_k(\bbc)$. Then, $\mathcal{B}^\prime\cap\mathcal{A}_1=\mbox{Alg}\{p_0,\ldots,p_{k-1}\}\otimes\bbc$ and $\mathcal{A}^\prime\cap\mathcal{A}_2=\mbox{Alg}\{I_k,C_k,\ldots,C_k^{k-1}\}\otimes\bbc$ are subalgebras of $M_k(\bbc)\otimes\mathscr{B}\left(\ell^2(\bbz^2)\right)$. The Fourier and the inverse Fourier transform are given by the following maps,
\begin{center}$
\mathcal{F}:\sum_{r=0}^{k-1}\alpha_rp_r\longmapsto\frac{1}{\sqrt{k}}\sum_{r=0}^{k-1}\alpha_rC_k^r\qquad\mbox{and}\qquad \mathcal{F}^{-1}:\sum_{r=0}^{k-1}\gamma_rC_k^r\longmapsto\sqrt{k}\,\sum_{r=0}^{k-1}\gamma_rp_r$
\end{center}
where $\,\alpha_r,\gamma_r\in\bbc$. Let us consider the following multiplication on $\bbc^k$,
\begin{eqnarray}\label{new mult}
(\alpha_0,\ldots,\alpha_{k-1})*(\beta_0,\ldots,\beta_{k-1}):=(\gamma_0,\ldots,\gamma_{k-1})
\end{eqnarray}
where $\gamma_j=\sum_{r=0}^{k-1}\alpha_r\beta_{k+j-r}$ for $0\leq j\leq k-1$, with the convention $\beta_{k+j}=\beta_j$ for all $j$. Then, $\mbox{Alg}\{I_k,C_k,\ldots,C_k^{k-1}\}\cong(\bbc^k,*)$ as unital algebras and the following map
\begin{center}$
\Phi:(\alpha_0,\ldots,\alpha_{k-1})\longmapsto\frac{1}{\sqrt{k}}\Big(\sum_{r=0}^{k-1}\alpha_r,\sum_{r=0}^{k-1}\omega^{r}\alpha_r,\sum_{r=0}^{k-1}\omega^{2r}\alpha_r,\ldots,\sum_{r=0}^{k-1}\omega^{(k-1)r}\alpha_r\Big)\,,$
\end{center}
where $\omega=e^{2\pi i/k}$ is a primitive $k$-th root of unity, implements a unital algebra isomorphism between $(\bbc^k,*)$ and $\bbc^k$ equipped with the standard algebra structure. It now follows that $\Phi\circ\mathcal{F}:\bbc^k\to\bbc^k$ is equal to the Fourier matrix. Moreover, the rotation maps $\rho_+,\rho_-$ are given by the following,
\begin{center}$
\rho_+:\sum_{r=0}^{k-1}\alpha_rp_r\longmapsto \sum_{r=1}^k\alpha_{k-r}p_r\quad\mbox{and}\quad\rho_-:\sum_{r=0}^{k-1}\gamma_rC_k^r\longmapsto \sum_{r=1}^k\gamma_{k-r}C_k^r\,,$
\end{center}
with the convention $p_k=p_0$. Therefore, as an element of $M_k(\bbc)$, both $\rho_+$ and $\,\Phi\circ\rho_-\circ\Phi^{-1}$ are equal to the $k\times k$ permutation matrix $E_{11}+\sum_{j=2}^kE_{j,k+2-j}$ in $M_k(\bbc)$. It turns out that the convolution product on $Alg\{p_0,\ldots,p_{k-1}\}\cong\bbc^k$ is the product $*$ defined in \Cref{new mult}, and that on $Alg\{I_k,C_k,\ldots,C_k^{k-1}\}\cong(\bbc^k,*)$ is the usual componentwise multiplication on $\bbc^k$.
\medskip

\subsection{Fourier transform for matrix algebras}

Let us consider the inclusion $\mathbb{C}\subset M_n(\mathbb{C})$. It is well known that the (standard normalized) trace-preserving conditional expectation is of index-finite type and it is the unique minimal one.  We shall use the notation $(\alpha_{ij})_{ij},\,1\leq i,j\leq n,$ to denote a matrix in $M_n(\bbc)$, and the elementary matrices will be denoted by $E_{ij}$. The unique normalized trace on $M_n(\bbc)$ is denoted by $\tr$. It is known that the basic construction for the unital inclusion $\mathbb{C}\subset M_n(\mathbb{C})$ is of the following form
\begin{center}$
\mathbb{C}\subset M_n(\mathbb{C})\subset^{\,e_1} M_n(\mathbb{C})\otimes M_n(\mathbb{C})\subset^{\,e_2} M_n(\mathbb{C})\otimes M_n(\mathbb{C})\otimes M_n(\mathbb{C})\subset\cdots\cdots$
\end{center}
with $e_1=\frac{1}{n}\sum_{i,j=1}^nE_{ij}\otimes E_{ij}$ and $e_2=\frac{1}{n}\sum_{i,j=1}^nE_{ij}\otimes E_{ij}\otimes I_n$ (see \cite{JS}). Thus, in accordance with the notations used in earlier sections, we have in this situation the inclusion $\mathcal{B}\subset\mathcal{A}\subset\mathcal{A}_1\subset\mathcal{A}_2$ where,
\[
\mathcal{B}=\mathbb{C}\otimes\mathbb{C}\otimes\mathbb{C}\quad,\quad\mathcal{A}=\mathbb{C}\otimes\mathbb{C}\otimes M_n(\mathbb{C})\,,
\]
\[
\mathcal{A}_1=\mathbb{C}\otimes M_n(\mathbb{C})\otimes M_n(\mathbb{C})\quad,\quad\mathcal{A}_2=M_n(\mathbb{C})\otimes M_n(\mathbb{C})\otimes M_n(\mathbb{C})\,.
\]
Clearly, $\mathcal{B}^\prime\cap\mathcal{A}_1=\mathcal{A}_1$ and $\mathcal{A}^\prime\cap\mathcal{A}_2=M_n(\mathbb{C})\otimes M_n(\mathbb{C})\otimes\mathbb{C}$. Then, it is clear that the conditional expectation $E^{\mathcal{A}_2}_{\mathcal{A}^\prime\cap\mathcal{A}_2}$ is given by $\text{id}\otimes\text{id}\otimes\tr$ and the conditional expectation $E^{\mathcal{A}_2}_{\mathcal{A}_1}$ is given by $\tr\otimes\text{id}\otimes\text{id}$. We shall use the standard convention $E_{(i,p)(j,q)}:=E_{ij}\otimes E_{pq}$ for the matrix units in $M_n(\mathbb{C})\otimes M_n(\mathbb{C})$ (Sec. $6$, Page $97$ in \cite{PP}).

\begin{ppsn}\label{to refer 1}
The Fourier and the inverse Fourier transform on $M_n(\mathbb{C})\otimes M_n(\mathbb{C})$ are given by the following,
\begin{center}
$\mathcal{F}:E_{k\ell}\otimes E_{pq}\longmapsto E_{\ell q}\otimes E_{kp}\,,$
\end{center}
\begin{center}
$\mathcal{F}^{-1}:E_{k\ell}\otimes E_{pq}\longmapsto E_{pk}\otimes E_{q\ell}\,.$
\end{center}
\end{ppsn}
We now find convolution on $M_n(\mathbb{C})\otimes M_n(\mathbb{C})$ defined by the formula $x* y:=\mathcal{F}^{-1}(\mathcal{F}(y)\mathcal{F}(x))$. Given two elements $A,D\in M_n(\mathbb{C})$, let $A\odot D$ be their Schur product. Since the projection $\frac{1}{n}J_n\in M_n(\mathbb{C})$, with $J_n=\sum_{i,j=1}^nE_{ij}$, is a minimal projection, we have $\frac{1}{n^2}J_n(A\odot D)J_n$ is a scalar multiple of $\frac{1}{n}J_n$. Denote this scalar by $\alpha_{A,D}$. That is,
\begin{eqnarray}\label{scalar multiple}
\big(\frac{1}{n}J_n\big)(A\odot D)\big(\frac{1}{n}J_n\big) &=& \alpha_{A,D}\big(\frac{1}{n}J_n\big)
\end{eqnarray}
with $\alpha_{A,D}\in\bbc$.

\begin{ppsn}\label{coproductmatrices}
The convolution on $M_n(\mathbb{C})\otimes M_n(\mathbb{C})$ implemented by the Fourier and the inverse Fourier transform is the given by the following,
\begin{center}$
(A\otimes B)*(C\otimes D)=n\,\alpha_{A,D}(C\otimes B)$
\end{center}
where $\alpha_{A,D}\in\mathbb{C}$ is as defined in Equation (\ref{scalar multiple}).
\end{ppsn}

\begin{ppsn}
The rotation maps $\rho_{+},\,\rho_{-}:M_n(\bbc)\otimes M_n(\bbc)\longrightarrow M_n(\bbc)\otimes M_n(\bbc)$ coincide, i,e. $\rho_{+}=\rho_{-}$, and they are given by the following,
\begin{center}$
E_{ij}\otimes E_{k\ell}\longmapsto E_{\ell k}\otimes E_{ji}\,.$
\end{center}
\end{ppsn}

\begin{rmrk}\label{rm}\rm
It is easy to see that $\rho_+$ is trace preserving, i,e., $\tr(\rho_+(x))=\tr(x)$ for all $x\in M_n(\bbc)\otimes M_n(\bbc)$. Therefore, the Young's inequality in \Cref{Younginequalitytheorem} becomes the following,
\begin{center}
$\|x* y\|_r\leq n^2\|x\|_p\|y\|_q$
\end{center}
for $x,y\in M_n(\bbc)\otimes M_n(\bbc)$, where $1\leq p,q,r\leq \infty$ and $\frac{1}{p}+\frac{1}{q}=\frac{1}{r}+1$.
\end{rmrk}
\bigskip

\section*{Acknowledgements}
We sincerely thank Zhengwei Liu for useful exchange and suggestions that eventually lead to a substantial improvement of the paper. K.C.B acknowledges the support of INSPIRE Faculty grant DST/INSPIRE/04/2019/002754. S.G acknowledges the support of SERB grant MTR/2021/000818. SM was affiliated to CMI when the majority of the work was carried out. She would like to thank CMI for the institute postdoctoral fellowship. 
\bigskip

%\bibliographystyle{amsplain}
%\bibliography{BGS_arxiv_version2}

\begin{thebibliography}{CS79}
	
	\bibitem{BakshiVedlattice}
	{\scshape Bakshi, Keshab Chandra; Gupta, Ved Prakash}. Lattice of intermediate subalgebras. {\em J. Lond. Math. Soc.} (2) \textbf{104} (2021), no.~5, 2082--2127, MR4368671, Zbl 07652699.
	
	\bibitem{Bisch}
	{\scshape Bisch, Dietmar}. A note on intermediate subfactors. {\em Pacific J. Math.} \textbf{163} (1994), no.~2, 201--216, MR1262294, Zbl 0814.46053.
	
	\bibitem{Bi94}
	{\scshape Bisch, Dietmar}. Bimodules, higher relative commutants and the fusion algebra associated to a subfactor. Operator algebras and their applications ({W}aterloo, {ON}, 1994/1995), Fields Inst. Commun., vol. \textbf{{13}},
	Amer. Math. Soc., Providence, RI, 1997, pp.~13--63, MR1424954, Zbl 0894.46046.
	
	\bibitem{BiJo2}
	{\scshape Bisch, Dietmar; Jones, V.F.R}. Singly generated planar algebras of small dimension. {\em Duke Math. J.} \textbf{101} (2000), no.~1, 41--75, MR1733737, Zbl 1075.46053.
	
	\bibitem{connesgeometry}
	{\scshape Connes, Alain}. Noncommutative geometry. Academic Press, Inc., San Diego, CA, 1994, MR1303779, Zbl 0818.46076.
	
	\bibitem{liuquantumsmooth}
	{\scshape Huang, Linzhe; Liu, Zhengwei; Wu, Jinsong}. Quantum smooth uncertainty principles for von {N}eumann bi-algebras. arXiv:2107.09057.
	
	\bibitem{liuquantumfrobenius}
	{\scshape Huang, Linzhe; Liu, Zhengwei; Wu, Jinsong}. Quantum convolution inequalities on {F}robenius von {N}eumann algebras. arXiv:2204.04401.
	
	\bibitem{Liunoncommutative}
	{\scshape Jiang, Chunlan; Liu, Zhengwei; Wu, Jinsong}. Noncommutative uncertainty principles. {\em J. Funct. Anal.} \textbf{270} (2016), no.~1, 264--311, MR3419762, Zbl 1352.46062.
	
	\bibitem{liuquantumgroup}
	{\scshape Jiang, Chunlan; Liu, Zhengwei; Wu, Jinsong}. Uncertainty principles for locally compact quantum groups. {\em J. Funct. Anal.} \textbf{274} (2018), no.~8, 2399--2445, MR3767437, Zbl 1403.43002.
	
	\bibitem{Jo}
	{\scshape Jones, V.F.R}. Index for subfactors. {\em Invent. Math.} \textbf{72} (1983), no.~1, 1--25, MR0696688, Zbl 0508.46040.
	
	\bibitem{Jo2}
	{\scshape Jones, V.F.R}. Planar algebras {I}. {\em New Zealand J. Math.} \textbf{52} (2021), 1--107, MR4374438, Zbl 1484.46067.
	
	\bibitem{JS}
	{\scshape Jones, V.F.R.; Sunder, V.S.} Introduction to subfactors. vol.~{\textbf{234}}, Cambridge University Press, 1997, MR1473221, Zbl 0903.46062. 
	
	\bibitem{kadisonschwarz}
	{\scshape Kadison, Richard V}. Non-commutative conditional expectations and their applications. Contemporary Mathematics \textbf{365} (2004), 143--180, MR2106820, Zbl 1080.46044.
	
	\bibitem{KajiwaraWatatani}
	{\scshape Kajiwara, Tsuyoshi; Watatani, Yasuo}. Jones index theory by {H}ilbert {$C^*$}-bimodules and {$K$}-theory. {\em Trans. Amer. Math. Soc.} \textbf{352} (2000), no.~8, 3429--3472, MR1624182, Zbl 0954.46034.
	
	\bibitem{kosaki}
	 {\scshape Kosaki, Hideki}. Applications of the complex interpolation method to a von {N}eumann algebra: noncommutative {$L^{p}$}-spaces. {\em J. Funct. Anal.} \textbf{56} (1984), no.~1, 29--78, MR0735704, Zbl 0604.46063.
	
	\bibitem{liu2016exchange}
	{\scshape Liu, Zhengwei}. Exchange relation planar algebras of small rank. {\em Trans. Amer. Math. Soc.} \textbf{368} (2016), no.~12, 8303--8348, MR3551573, Zbl 1365.46053. 
	
	\bibitem{liukacalgebra}
	{\scshape Liu, Zhengwei; Wu, Jinsong}. Extremal pairs of {Y}oung’s inequality for {K}ac algebras. {\em Pacific Journal of Mathematics} \textbf{295} (2018), no.~1, 103--121, MR3778328, Zbl 1395.46056. 
	
	\bibitem{PP}
	{\scshape Pimsner, Mihai; Popa, Sorin}. Entropy and index for subfactors. {\em Ann. Sci. \'{E}cole Norm. Sup.} (4), vol.~{\textbf{19}}, 1986, pp.~57--106, MR0860811, Zbl 0646.46057. 
	
	\bibitem{P}
	{\scshape Popa, Sorin}. An axiomatization of the lattice of higher relative commutants of a subfactor. {\em Invent. Math.} \textbf{120} (1995), no.~3, 427--445, MR1334479, Zbl 0831.46069.
	
	\bibitem{stratila}
	{\scshape Stratila, S.V.; Zsido, L}. Lectures on von {N}eumann algebras. {C}ambridge {IIS}c {S}eries, (2019), Zbl 1418.46002.
	
	\bibitem{Watataniindex}
	 {\scshape Watatani, Yasuo}. Index for {$C^*$}-subalgebras. {\em Mem. Amer. Math. Soc.} \textbf{83} (1990), no.~424, vi+117, MR0996807, Zbl 0697.46024.
	
	\bibitem{Xunotes}
	{\scshape Xu, Q}. Operator spaces and noncommutative ${L}_p$: the part on noncommutative ${L}_p$-spaces, in: Lectures in the Summer School on ``Banach Spaces and Operator Spaces". Tianjin: Nankai University, \textbf{66} (2007).
\end{thebibliography}
%\end{document}

\bigskip

\bigskip

\noindent{\sc Keshab Chandra Bakshi} (\texttt{keshab@iitk.ac.in, bakshi209@gmail.com})\\
         {\footnotesize Department of Mathematics and Statistics,\\
         Indian Institute of Technology, Kanpur,\\
         Uttar Pradesh 208016, India.}
\bigskip

\noindent{\sc Satyajit Guin} (\texttt{sguin@iitk.ac.in})\\
         {\footnotesize Department of Mathematics and Statistics,\\
         Indian Institute of Technology, Kanpur,\\
         Uttar Pradesh 208016, India.}
\bigskip

\noindent{\sc Sruthymurali} (\texttt{sruthy92smk@gmail.com, sruthymurali\_pd@isibang.ac.in})\\
         {\footnotesize Indian Statistical Institute,\\ 8th Mile, Mysore Rd, RVCE Post,\\  Bengaluru, Karnataka 560059, India.}

\end{document}